# Modular Analysis of Almost Block Diagonal Systems of Equations

**Tarek M. A. El-Mistikawy**

Department of Engineering Mathematics and Physics, Cairo University, Giza 12211, Egypt

## Abstract

Almost block diagonal linear systems of equations can be exemplified by two modules. This makes it possible to construct all sequential forms of band and/or block elimination methods, six old and fourteen new. It allows easy assessment of the methods on the basis of their operation counts, storage needs, and admissibility of partial pivoting. It unveils a robust partial pivoting strategy- *local pivoting*. Extension of modular analysis to bordered systems is also included.

## 1. INTRODUCTION

Systems of equations with almost block diagonal (ABD) matrix of coefficients are frequently encountered in numerical solutions of sets of ordinary or partial differential equations. They occur, for example, when Keller's box scheme is applied: to parabolic equations, and to elliptic equations with an alternating-direction-implicit approach. Several other situations in which ABD systems occur are described by Amodio et al. [1].

Sequential solution methods for ABD systems were developed, improved, and thoroughly studied to the extent that it was concluded in Ref. [1] that they needed little further study.

Traditionally, sequential solution methods of ABD systems performed $\mathbb{LU}$ decompositions of the matrix of coefficients $\mathbb{G}$ through either band (scalar) elimination or block tridiagonal elimination. The famous COLROW algorithm [2], which is highly regarded for its performance, was incorporated in several applications [3,4,5,6,7]. It utilizes Lam's alternating column/row pivoting [8] and Varah's correspondingly alternating column/row scalar elimination [9]. The efficient block tridiagonal methods included Keller's Block Tridiagonal Row elimination method [10, §5, case i], and El-Mistikawy's Block Tridiagonal Column elimination method [11]. Both methods could apply a suitable form of Keller's mixed pivoting strategy [10], which is more expensive than Lam's.

The present paper is intended to explore other variants of the $\mathbb{LU}$ decomposition of $\mathbb{G}$. It does not follow the traditional approaches of treating the matrix of coefficients as a banded matrix, or casting it into a block tridiagonal form. It, rather, adopts a new approach, *modular analysis*, which offers a simple and unified way of expressing and assessing solution methods for ABD systems.

The matrix of coefficients $\mathbb{G}$ (or, more specifically, its significant part containing the non-zero blocks) is disassembled into an ordered set of modules. (In fact, two different sets of modules are identified.) Each module $\Gamma$ is an entity that has the following components: a stem $\sigma$, two fins $\psi$ and $\phi$, a head $\eta$, and a tail $\tau$. By arranging the modules in such a way that the head of a module is added to the tail of the next, the significant part of $\mathbb{G}$ can be re-assembled. The module exemplifies the matrix, but is much easier to analyze.

All possible methods of $\mathbb{LU}$ decomposition of $\mathbb{G}$ can be formulated as decompositions of Γ. Fourteen new methods are thus discovered. The operation counts and storage needs are easily estimated, revealing the method with the best performance on each account.

The validity and stability of the elimination methods are of primary concern to both numerical analysts and algorithm users. Validity means that division by a zero is never encountered, whereas stability guards against round-off-error growth. To insure validity and achieve stability, pivoting is called for [12]. Full pivoting is computationally expensive. It requires full two dimensional search for the pivots. Moreover, it destroys the banded form of the matrix of coefficients. Partial pivoting strategies, though potentially less stable, are considerably less expensive. Uni-directional (row or column) pivoting makes a minor change to the form of $\mathbb{G}$ by introducing few extraneous elements. Lam's alternating pivoting [8], which involves alternating sequences of row pivoting and column pivoting, maintains the form of $\mathbb{G}$. When $\mathbb{G}$ is nonsingular, Lam's pivoting guarantees validity, and if followed by correspondingly alternating elimination it produces multipliers that are bounded by unity; thus enhancing stability. This approach was proposed by Varah [9] in his $\mathbb{L}\tilde{\mathbb{G}}\,\mathbb{U}$ decomposition method. It was developed afterwards into a more efficient $\mathbb{LU}$ version that was adopted by the well-known COLROW solver [2]. The present approach of modular analysis shows that Lam's pivoting (with Varah's arrangement) applies to five other elimination methods. A more robust, though more expensive strategy- *Local Pivoting*, is hereinafter introduced. It performs full pivoting over the same segments of $\mathbb{G}$ (or Γ) to which Lam's pivoting is applied.

## 2. PROBLEM DESCRIPTION

Consider the almost block diagonal system of equations $\mathbb{G}\mathsf{z} = \mathsf{g}$ whose augmented matrix of coefficients $\mathbb{G}^+ = [\mathbb{G} \;\vdots\; \mathsf{g}]$ has the form shown in Fig. 1, where the blocks with leading character *A*, *B*, *C*, and *D* have *m*, *n*, *m*, and *n* columns, respectively. The blocks with leading character *g* have *r* columns indicating as many right-hand sides. The trailing character *m* or *n* (and subsequently *p*=*m*+*n*) indicates the number of rows of a block; or the order of a square matrix such as an identity matrix *I*, a zero matrix *0*, and a lower *L* or an upper *U* triangular matrix. Blanks indicate zero blocks.

$$\mathbb{G}^+ = \left[\begin{array}{cccccccccc:c}
Cm_1 & Dm_1 & & & & & & & & & gm_1 \\
An_1 & Bn_1 & Cn_2 & Dn_2 & & & & & & & gn_1 \\
Am_1 & Bm_1 & Cm_2 & Dm_2 & & & & & & & gm_2 \\
 & \cdots & \cdots & \cdots & \cdots & & & & & & \cdots \\
 & \vdots & An_{j-1} & Bn_{j-1} & Cn_j & Dn_j & & & & & gn_{j-1} \\
 & \vdots & Am_{j-1} & Bm_{j-1} & Cm_j & Dm_j & & & & & gm_j \\
 & & & \cdots & \cdots & \cdots & \cdots & & & & \cdots \\
 & & & \vdots & An_{J-1} & Bn_{J-1} & Cn_J & Dn_J & & & gn_{J-1} \\
 & & & \vdots & Am_{J-1} & Bm_{J-1} & Cm_J & Dm_J & & & gm_J \\
 & & & & & & An_J & Bn_J & & & gn_J
\end{array}\right]$$

Fig. 1: The augmented matrix of coefficients

The matrix of unknowns $\mathsf{z}$ is, likewise, written as

$$\mathsf{z} = \begin{bmatrix} zm_1{}^t & zn_1{}^t & zm_2{}^t & zn_2{}^t & \cdots & zm_j{}^t & zn_j{}^t & \cdots & zm_J{}^t & zn_J{}^t \end{bmatrix}^t$$

where the superscript $t$ denotes the transpose.

Such a system of equations results, for example, in the finite difference solution of $p$ first order ODEs on a grid of $J$ points; with $m$ conditions at one side to be marked with $j=1$, and $n$ conditions at the other side that is to be marked with $j=J$. Then, each column of the submatrix $\begin{bmatrix} zm_j{}^t & zn_j{}^t \end{bmatrix}^t$ contains the $p$ unknowns of the $j^{\text{th}}$ grid point corresponding to a right-hand side.

## 3. MODULAR ANALYSIS

The description of the decomposition methods for the augmented matrix of coefficients $\mathbb{G}^+$ can be made easy and concise through the introduction of modules of $\mathbb{G}^+$. Two different modules can be identified as described below.

### 3.1. The Aligned Module (A-Module)

This module, for $j$=1→$J$-1, takes the form

$$\Gamma_j^{+A} \equiv \left[\Gamma_j^A \;\vdots\; \gamma_j^A\right] = \left[\begin{array}{cc:cc:c} Cm_j^{\#} & Dm_j^{\#} & & & gm_j^{\#} \\ An_j & Bn_j & Cn_{j+1} & Dn_{j+1} & gn_j \\ \hdashline Am_j & Bm_j & Cm_{j+1}^{\Rightarrow} & Dm_{j+1}^{\Rightarrow} & gm_{j+1}^{\Rightarrow} \end{array}\right]$$

where, as a rule, the dotted line defines the partitioning to left and right entities. The dashed lines define the partitioning

$$\Gamma_j^{+A} \equiv \left[\begin{array}{c:c} \sigma_j^A & \phi_j^{+A} \\ \hdashline \psi_j^A & \eta_j^{+A} \end{array}\right]$$

of $\Gamma_j^{+A}$ to the following components.

The stem

$$\sigma_j^A = \begin{bmatrix} Cm_j^{\#} & Dm_j^{\#} \\ An_j & Bn_j \end{bmatrix}$$

The head

$$\eta_j^{+A} \equiv \left[\eta_j^A \;\vdots\; \theta_j^A\right] = \left[Cm_{j+1}^{\Rightarrow} \quad Dm_{j+1}^{\Rightarrow} \;\vdots\; gm_{j+1}^{\Rightarrow}\right]$$

The fins

$$\psi_j^A = \begin{bmatrix} Am_j & Bm_j \end{bmatrix}$$

and

$$\phi_j^{+A} \equiv \left[\phi_j^A \;\vdots\; \rho_j^A\right] = \begin{bmatrix} & & \vdots & gm_j^{\#} \\ Cn_{j+1} & Dn_{j+1} & \vdots & gn_j \end{bmatrix}$$

The module has the tail

$$\tau_j^{+A} \equiv \left[\tau_j^A \;\vdots\; \lambda_j^A\right] = \left[Cm_j^{\#} \quad Dm_j^{\#} \;\vdots\; gm_j^{\#}\right]$$

which is defined through the head-tail relation

$$\eta_{j-1}^{+A} + \tau_j^{+A} = \left[Cm_j \quad Dm_j \;\vdots\; gm_j\right]$$

This makes it possible to construct the significant part of $\mathbb{G}^+$ by arranging the set of modules in such a way that the tail of $\Gamma_j^{+A}$ adds to the head of $\Gamma_{j-1}^{+A}$. Minor adjustments need only to be invoked at both ends of $\mathbb{G}^+$. Specifically, we define the truncated modules $\Gamma_0^{+A} \equiv \eta_0^{+A} = 0$ and $\Gamma_J^{+A} \equiv \left[\sigma_J^A \;\vdots\; \rho_J^A\right]$.

(For convenience, we shall occasionally drop the subscript and/or superscript identifying a module and its components, as well as the subscript identifying its blocks. This should not cause any confusion.)

The head of the module $\Gamma^+$ is yet to be defined. It is taken to be related to the other components of $\Gamma^+$ by

$$\eta^+ = \psi\sigma^{-1}\phi^+ \tag{3.1}$$

in order to allow for decompositions of $\Gamma^+$ having the form

$$\Gamma^+ = \mu\nu^+ \equiv \begin{bmatrix} \mathrm{M} \\ \text{- -} \\ \Psi \end{bmatrix} \left[\mathrm{N} \;\vdots\; \Phi^+\right] \tag{3.2}$$

The generic relations $\sigma$=MN, $\psi$=ΨN, $\phi^+ = \mathrm{M}\Phi^+$, and $\eta^+ = \Psi\Phi^+$ then hold; leading to $\eta^+ = \Psi\Phi^+ = (\psi\mathrm{N}^{-1})(\mathrm{M}^{-1}\phi^+) = \psi(\mathrm{N}^{-1}\mathrm{M}^{-1})\phi^+ = \psi\sigma^{-1}\phi^+$ as defined in (3.1). Further, writing $\Phi^+ \equiv \left[\Phi \;\vdots\; \mathrm{Z}\right]$ gives $\phi$=MΦ, $\rho$=MZ, $\eta$=ΨΦ, and $\theta$=ΨZ.

We are now ready to present several old as well as new (marked with a square ■) elimination methods in terms of decompositions of the left module $\Gamma_j^A$. Several inflections of the blocks of $\mathbb{G}$ are involved and are defined in Appendix A. The sequence in which the blocks are manipulated for: decomposing the stem, processing the fins, and handling the head (evaluating the head and applying the head-tail relation to determine the tail of the succeeding module $\Gamma_{j+1}^A$), is mentioned along with the equations (from Appendix A) involved. The correctness of the decompositions may be checked by carrying out the matrix multiplications, using the equalities of Appendix A, and comparing with the un-decomposed form of the module.

### 3.1.1. Scalar/Scalar Elimination Methods

These methods perform scalar decomposition of both pivotal blocks $Cm^{\#}$ and $Bn^{\#}$. The triangular matrices $L$ and $U$ appear explicitly, marked if unit diagonal with a circumflex $\hat{}$ . Four methods are involved.

**Scalar Row/Scalar Row (SRSR) Elimination**
(Case iii of [10, §5])

$$\Gamma^A = \begin{bmatrix} \hat{L}m & \\ An' & \hat{L}n \\ \hdashline Am' & Bm'' \end{bmatrix} \left[\begin{array}{cc:cc} Um & Dm'' & & \\ & Un & Cn' & Dn' \end{array}\right]$$

**Scalar Row/Scalar Column (SRSC) Elimination■**

$$\Gamma^A = \begin{bmatrix} \hat{L}m & \\ An' & Ln \\ \hdashline Am' & Bm'' \end{bmatrix} \left[\begin{array}{cc:cc} Um & Dm'' & & \\ & \hat{U}n & Cn' & Dn' \end{array}\right]$$

**Scalar Column/Scalar Row (SCSR) Elimination**
(Method implemented by the COLROW algorithm)

$$\Gamma^A = \begin{bmatrix} Lm & \\ An' & \hat{L}n \\ \hdashline Am' & Bm'' \end{bmatrix} \left[\begin{array}{cc:cc} \hat{U}m & Dm'' & & \\ & Un & Cn' & Dn' \end{array}\right]$$

**Scalar Column/Scalar Column (SCSC) Elimination**

(Case iv of [10, §5])

$$\Gamma^A = \begin{bmatrix} Lm & \\ An' & Ln \\ \hdashline Am' & Bm'' \end{bmatrix} \left[\begin{array}{cc:cc} \hat{U}m & Dm'' & & \\ & \hat{U}n & Cn' & Dn' \end{array}\right]$$

All four scalar elimination methods apply the following sequence of manipulations.

Stem: $LmUm$(A6), $Dm''$(A19), $An'$(A9), $Bn^{\#}$ (A2b), $LnUn$(A5)
Fins: $Am'$(A7), $Bm^{\#}$ (A1b), $Bm''$(A17), $Cn'$(A11), $Dn'$(A13)
Head: $Cm^{\#}$ (A3b), $Dm^{\#}$ (A4b)

The difference is in how the decompositions of the pivotal blocks $Cm^{\#}$ and $Bn^{\#}$ are carried out, through row ( $\hat{L}U$ ) or column ( $L\hat{U}$ ) manipulation.

### 3.1.2. Scalar/Block Elimination Methods

These methods perform scalar decomposition of $Cm^{\#}$ and block decomposition of $Bn^{\#}$. It involves four methods.

**Scalar Row/Block Row (SRBR) Elimination■**

$$\Gamma^{A} = \left[\begin{array}{cc} \hat{L}m & \\ An' & In \\ \hdashline Am' & Bm^{*} \end{array}\right] \left[\begin{array}{cc:cc} Um & Dm'' & & \\ & \breve{B}n^{\#} & Cn & Dn \end{array}\right]$$

**Scalar Column/Block Row (SCBR) Elimination■**

$$\Gamma^{A} = \left[\begin{array}{cc} Lm & \\ An' & In \\ \hdashline Am' & Bm^{*} \end{array}\right] \left[\begin{array}{cc:cc} \hat{U}m & Dm^{\#} & & \\ & \breve{B}n^{\#} & Cn & Dn \end{array}\right]$$

Both SRBR and SCBR methods apply the following sequence of manipulations.

Stem: $\breve{C}m^{\#}$ (A6), $An'$(A9), $Dm''$(A19), $Bn^{\#}$ (A2b), $\breve{B}n^{\#}$ (A5)
Fins: $Am'$(A7) , $Bm^{\#}$ (A1b), $Bm''$(A17), $Bm^{*}$(A18)
Head: $Cm^{\#}$ (A3a), $Dm^{\#}$ (A4a)

The difference is in how the decompositions of the pivotal block $Cm^{\#}$ is carried out, through row ( $\hat{L}U$ ) or column ( $L\hat{U}$ ) manipulation.

**Scalar Row/Block Column (SRBC) Elimination■**

$$\Gamma^{A} = \left[\begin{array}{cc} \hat{L}m & \\ An' & \breve{B}n^{\#} \\ \hdashline Am' & Bm^{\#} \end{array}\right] \left[\begin{array}{cc:cc} Um & Dm'' & & \\ & In & Cn^{\circ} & Dn^{\circ} \end{array}\right]$$

**Scalar Column/Block Column (SCBC) Elimination■**

$$\Gamma^{A} = \left[\begin{array}{cc} Lm & \\ An' & \breve{B}n^{\#} \\ \hdashline Am' & Bm^{\#} \end{array}\right] \left[\begin{array}{cc:cc} \hat{U}m & Dm^{\#} & & \\ & In & Cn^{\circ} & Dn^{\circ} \end{array}\right]$$

Both SRBC and SCBC methods apply the following sequence of manipulations.

Stem: $\breve{C}m^{\#}$ (A6), $An'$(A9), $Dm''$(A19), $Bn^{\#}$ (A2b), $\breve{B}n^{\#}$ (A5)
Fins: $Am'$(A7), $Bm^{\#}$ (A1b), $Cn'$(A11), $Dn'$(A13), $Cn^{\circ}$(A12), $Dn^{\circ}$(A14)
Head: $Cm^{\#}$ (A3c), $Dm^{\#}$ (A4c)

The difference is in how the decompositions of the pivotal block $Cm^{\#}$ is carried out, through row ( $\hat{L}U$ ) or column ( $L\hat{U}$ ) manipulation.

### 3.1.3. Block/Scalar Elimination Methods

These methods perform block decomposition of $Cm^{\#}$ and scalar decomposition of $Bn^{\#}$. It involves four methods.

**Block Row/Scalar Row (BRSR) Elimination**■

$$\Gamma^{A} = \begin{bmatrix} Im & \\ An^{\circ} & \hat{L}n \\ \hdashline Am^{\circ} & Bm'' \end{bmatrix} \begin{bmatrix} \breve{C}m^{\#} & Dm^{\#} & \vdots & & \\ & Un & \vdots & Cn' & Dn' \end{bmatrix}$$

**Block Row/Scalar Column (BRSC) Elimination**■

$$\Gamma^{A} = \begin{bmatrix} Im & \\ An^{\circ} & Ln \\ \hdashline Am^{\circ} & Bm'' \end{bmatrix} \begin{bmatrix} \breve{C}m^{\#} & Dm^{\#} & \vdots & & \\ & \hat{U}n & \vdots & Cn' & Dn' \end{bmatrix}$$

Both BRSR and BRSC methods apply the following sequence of manipulations.

Stem: $\breve{C}m^{\#}$ (A6), $An'$(A9), $An^{\circ}$(A10), $Bn^{\#}$ (A2c), $\breve{B}n^{\#}$ (A5)
Fins: $Am'$(A7), $Am^{\circ}$(A8), $Bm^{\#}$ (A1c), $Bm''$(A17)
Head: $Cm^{\#}$ (A3b), $Dm^{\#}$ (A4b)

The difference is in how the decompositions of the pivotal block $Bn^{\#}$ is carried out, through row ( $\hat{L}U$ ) or column ( $L\hat{U}$ ) manipulation.

**Block Column/Scalar Row (BCSR) Elimination**■

$$\Gamma^{A} = \begin{bmatrix} \breve{C}m^{\#} & \\ An & \hat{L}n \\ \hdashline Am & Bm'' \end{bmatrix} \begin{bmatrix} Im & Dm^{*} & \vdots & & \\ & Un & \vdots & Cn' & Dn' \end{bmatrix}$$

**Block Column/Scalar Column (BCSC) Elimination**■

$$\Gamma^{A} = \begin{bmatrix} \breve{C}m^{\#} & \\ An & Ln \\ \hdashline Am & Bm'' \end{bmatrix} \begin{bmatrix} Im & Dm^{*} & \vdots & & \\ & \hat{U}n & \vdots & Cn' & Dn' \end{bmatrix}$$

Both BCSR and BCSC methods apply the following sequence of manipulations.

Stem: $\breve{C}m^{\#}$ (A6), $Dm''$(A19), $Dm^{*}$(A20), $Bn^{\#}$ (A2a), $\breve{B}n^{\#}$ (A5)
Fins: $Bm^{\#}$ (A1a), $Bm''$(A17), $Cn'$(A11), $Dn'$(A13)
Head: $Cm^{\#}$ (A3b), $Dm^{\#}$ (A4b)

The difference is in how the decompositions of the pivotal block $Bn^{\#}$ is carried out, through row ( $\hat{L}U$ ) or column ( $L\hat{U}$ ) manipulation.

### 3.1.4. Block/Block Elimination Methods

These methods perform block decomposition of both pivotal blocks $Cm^{\#}$ and $Bn^{\#}$, in which the decomposed pivotal blocks $\breve{C}m^{\#}$ ($\equiv LmUm$) and $\breve{B}n^{\#}$ ($\equiv LnUn$) appear. It involves four methods.

**Block Row/Block Row (BRBR) Elimination**■

$$\Gamma^{A} = \left[\begin{array}{cc} Im & \\ An^{\circ} & In \\ \hdashline Am^{\circ} & Bm^{*} \end{array}\right] \left[\begin{array}{cc:cc} \breve{C}m^{\#} & Dm^{\#} & & \\ & \breve{B}n^{\#} & Cn & Dn \end{array}\right]$$

Stem: $\breve{C}m^{\#}$ (A6), $An'$(A9), $An^{\circ}$(A10), $Bn^{\#}$ (A2c), $\breve{B}n^{\#}$ (A5)
Fins: $Am'$(A7), $Am^{\circ}$(A8), $Bm^{\#}$ (A1c), $Bm''$(A17), $Bm^{*}$(A18)
Head: $Cm^{\#}$ (A3a), $Dm^{\#}$ (A4a)

**Block Row/Block Column (BRBC) Elimination**■

$$\Gamma^{A} = \left[\begin{array}{cc} Im & \\ An^{\circ} & \breve{B}n^{\#} \\ \hdashline Am^{\circ} & Bm^{\#} \end{array}\right] \left[\begin{array}{cc:cc} \breve{C}m^{\#} & Dm^{\#} & & \\ & In & Cn^{\circ} & Dn^{\circ} \end{array}\right]$$

Stem: $\breve{C}m^{\#}$ (A6), $An'$(A9), $An^{\circ}$(A10), $Bn^{\#}$ (A2c), $\breve{B}n^{\#}$ (A5)
Fins: $Am'$(A7), $Am^{\circ}$(A8), $Bm^{\#}$ (A1c), $Cn'$(A11), $Dn'$(A13), $Cn^{\circ}$(A12), $Dn^{\circ}$(A14)
Head: $Cm^{\#}$ (A3c), $Dm^{\#}$ (A4c)

**Block Column/Block Row (BCBR) Elimination**■

$$\Gamma^{A} = \left[\begin{array}{cc} \breve{C}m^{\#} & \\ An & In \\ \hdashline Am & Bm^{*} \end{array}\right] \left[\begin{array}{cc:cc} Im & Dm^{*} & & \\ & \breve{B}n^{\#} & Cn & Dn \end{array}\right]$$

Stem: $\breve{C}m^{\#}$ (A6), $Dm''$(A19), $Dm^{*}$(A20), $Bn^{\#}$ (A2a), $\breve{B}n^{\#}$ (A5)
Fins: $Bm^{\#}$ (A1a), $Bm''$(A17), $Bm^{*}$(A18)
Head: $Cm^{\#}$ (A3a), $Dm^{\#}$ (A4a)

**Block Column/Block Column (BCBC) Elimination**■

$$\Gamma^{A} = \left[\begin{array}{cc} \breve{C}m^{\#} & \\ An & \breve{B}n^{\#} \\ \hdashline Am & Bm^{\#} \end{array}\right] \left[\begin{array}{cc:cc} Im & Dm^{*} & & \\ & In & Cn^{\circ} & Dn^{\circ} \end{array}\right]$$

Stem: $\breve{C}m^{\#}$ (A6), $Dm''$(A19), $Dm^{*}$(A20), $Bn^{\#}$ (A2a), $\breve{B}n^{\#}$ (A5)
Fins: $Bm^{\#}$ (A1a), $Cn'$(A11), $Dn'$(A13), $Cn^{\circ}$(A12), $Dn^{\circ}$(A14)
Head: $Cm^{\#}$ (A3c), $Dm^{\#}$ (A4c)

### 3.1.5. Aligned Block-Tridiagonal Elimination Methods

These methods perform the identity decomposition $\sigma = Ip\breve{\sigma}$ or $\sigma = \breve{\sigma}Ip$, where the decomposed stem $\breve{\sigma}$ stands for any of the nonidentity (scalar or block) decompositions given above. This is an arbitrariness that can be put to advantage, as shown below.

**Aligned Block-Tridiagonal Row (ABTR) Elimination**
**(**Case i of [10, §5]**)**

$$\Gamma^A = \begin{bmatrix} Ip \\ \hdashline Am^* \quad Bm^* \end{bmatrix} \begin{bmatrix} \breve{\sigma}^A & \vdots & Im & \\ & \vdots & Cn & Dn \end{bmatrix}$$

Traditionally ABTR elimination decomposes the stem $\sigma^A$ through scalar row manipulation. However, as indicated in §5 on pivoting, it would prove profitable to decompose the stem through an alternating manipulation method as in !C?R elimination, where ! and ? stand for S or B. Using SCSR elimination, the following sequence of manipulations applies.

Stem: $Lm\hat{U}m$ (A6), $Dm''$(A19), $An'$(A9), $Bn^{\#}$ (A2b), $\hat{L}\,nUn$(A5)
Fins: $Am'$(A7), $Bm^{\#}$ (A1b), $Bm''$(A17), $Bm^*$(A18), $Am''$(A15), $Am^*$(A16)
Head: $Cm^{\#}$ (A3a), $Dm^{\#}$ (A4a)

**Aligned Block-Tridiagonal Column (ABTC) Elimination**
**(**Case ii of [10, §5]**)**

$$\Gamma^A = \begin{bmatrix} \breve{\sigma}^A \\ \hdashline Am \quad Bm \end{bmatrix} \begin{bmatrix} Ip & \vdots & Cm^e & Dm^e \\ & \vdots & Cn^\circ & Dn^\circ \end{bmatrix}$$

The method introduces the two extraneous blocks $Cm^e$ and $Dm^e$. Traditionally ABTC elimination (now, ABTCt) decomposes the stem $\sigma^A$ through scalar row manipulation as in SRSR elimination.

Stem: $\hat{L}\,mUm$(A6), $Dm''$(A19), $An'$(A9), $Bn^{\#}$ (A2b), $\hat{L}\,nUn$(A5)
Fins: $Cn'$(A11), $Dn'$(A13), $Cn^\circ$(A12), $Dn^\circ$(A14), $Cm^e$ (A25b), $Dm^e$ (A26b)
Head: $Cm^{\#}$ (A29), $Dm^{\#}$ (A30)

A more efficient sequence (defining ABTCe) decomposes the stem through block column manipulation as in BCBC elimination.

Stem: $\breve{C}m^{\#}$ (A6), $Dm''$(A19), $Dm^*$(A20), $Bn^{\#}$ (A2a), $\breve{B}n^{\#}$ (A5)
Fins: $Cn'$(A11), $Dn'$(A13), $Cn^\circ$(A12), $Dn^\circ$(A14), $Cm^e$ (A25a), $Dm^e$ (A26a)
Head: $Bm^{\#}$ (A1a), $Cm^{\#}$ (A3c), $Dm^{\#}$ (A4c)

Note that the use of (A3c) and (A4c) for evaluating the head is a result of substituting expressions (A25a) and (A26a) for $Cm^e$ and $Dm^e$ into (A29) and (A30), respectively.

### 3.2. The Displaced Module (D-Module)

The second module of $\mathbb{G}^+$ takes, for $j$=1⟶$J$-2, the form

$$\Gamma_j^{+D} \equiv \left[\Gamma_j^D \;\vdots\; \gamma_j^D\right] = \left[\begin{array}{cc:c:c} Bn_j^\# & Cn_{j+1} & Dn_{j+1} & gn_j^\# \\ Bm_j^\# & Cm_{j+1} & Dm_{j+1} & gm_{j+1}^\# \\ \hdashline & An_{j+1} & Bn_{j+1}^\Rightarrow & gn_{j+1}^\Rightarrow \\ & Am_{j+1} & Bm_{j+1}^\Rightarrow & gm_{j+2}^\Rightarrow \end{array}\right]$$

with the partitioning

$$\Gamma_j^{+D} \equiv \left[\begin{array}{c:c} \sigma_j^D & \phi_j^{+D} \\ \hdashline \psi_j^D & \eta_j^{+D} \end{array}\right]$$

which defines the following components of $\Gamma_j^{+D}$.

The stem

$$\sigma_j^D = \begin{bmatrix} Bn_j^\# & Cn_{j+1} \\ Bm_j^\# & Cm_{j+1} \end{bmatrix}$$

The head

$$\eta_j^{+D} \equiv \left[\eta_j^D \;\vdots\; \theta_j^D\right] = \left[\begin{array}{c:c} Bn_{j+1}^\Rightarrow & gn_{j+1}^\Rightarrow \\ Bm_{j+1}^\Rightarrow & gm_{j+2}^\Rightarrow \end{array}\right]$$

The fins

$$\psi_j^D = \begin{bmatrix} 0n & An_{j+1} \\ & Am_{j+1} \end{bmatrix}$$

and

$$\phi_j^{+D} \equiv \left[\phi_j^D \;\vdots\; \rho_j^D\right] = \left[\begin{array}{c:c} Dn_{j+1} & gn_j^\# \\ Dm_{j+1} & gm_{j+1}^\# \end{array}\right]$$

The tail

$$\tau_j^{+D} \equiv \left[\tau_j^D \;\vdots\; \lambda_j^D\right] = \left[\begin{array}{c:c} Bn_j^\# & gn_j^\# \\ Bm_j^\# & gm_{j+1}^\# \end{array}\right]$$

The head-tail relation for $\Gamma_j^{+D}$ takes the form

$$\eta_{j-1}^{+D} + \tau_j^{+D} = \left[\begin{array}{c:c} Bn_j & gn_j \\ Bm_j & gm_{j+1} \end{array}\right]$$

so that, with minor end adjustments involving the use of the truncated modules

$$\Gamma_0^{+D}=\left[\begin{array}{c:c:c} Cm_1 & Dm_1 & gm_1 \\ \hdashline An_1 & Bn_1^{\Rightarrow} & gn_1^{\Rightarrow} \\ Am_1 & Bm_1^{\Rightarrow} & gm_2^{\Rightarrow} \end{array}\right],$$

$$\Gamma_{J-1}^{+D}=\left[\begin{array}{cc:c:c} Bn_{J-1}^{\#} & Cn_J & Dn_J & gn_{J-1}^{\#} \\ Bm_{J-1}^{\#} & Cm_J & Dm_J & gm_J^{\#} \\ \hdashline & An_J & Bn_J^{\Rightarrow} & gn_J^{\Rightarrow} \end{array}\right], \text{ and}$$

$$\Gamma_J^{+D}=\left[\begin{array}{c:c} Bn_J^{\#} & gn_J^{\#} \end{array}\right]$$

the significant part of $\mathbb{G}^{+}$ can be constructed by arranging the modules in such a way that the tail of $\Gamma_j^{+D}$ adds to the head of $\Gamma_{j-1}^{+D}$.

The head of $\Gamma^{+D}$ is defined in accordance with relation (3.1) to allow for decompositions of the form (3.2). Such decompositions of $\Gamma^{D}$ lead to all but the aligned block-tridiagonal elimination methods given above; e.g.,

**Scalar Column/Scalar Column (SCSC) Elimination**

$$\Gamma^{D}=\left[\begin{array}{cc} Ln & \\ Bm'' & Lm \\ \hdashline & An' \\ & Am' \end{array}\right]\left[\begin{array}{cc:c} \hat{U}n & Cn' & Dn' \\ & \hat{U}m & Dm'' \end{array}\right]$$

Stem: $Ln\hat{U}n$ (A5), $Cn'$(A11), $Bm''$(A17), $Cm^{\#}$ (A3b), $Lm\hat{U}m$ (A6)
Fins: $Dn'$(A13), $Dm^{\#}$ (A4b), $Dm''$(A19), $An'$(A9), $Am'$(A7)
Head: $Bn^{\#}$ (A2b), $Bm^{\#}$ (A1b)

**Block Row/Block Row (BRBR) Elimination**■

$$\Gamma^{D}=\left[\begin{array}{cc} In & \\ Bm^{*} & Im \\ \hdashline & An^{\circ} \\ & Am^{\circ} \end{array}\right]\left[\begin{array}{cc:c} \breve{B}n^{\#} & Cn & Dn \\ & \breve{C}m^{\#} & Dm^{\#} \end{array}\right]$$

Stem: $\breve{B}n^{\#}$ (A5), $Bm''$(A17), $Bm^{*}$(A18), $Cm^{\#}$ (A3a), $\breve{C}m^{\#}$ (A6)
Fins: $Dm^{\#}$ (A4a), $An'$(A9), $Am'$(A7), $An^{\circ}$(A10), $Am^{\circ}$(A8)
Head: $Bn^{\#}$ (A2c), $Bm^{\#}$ (A1c)

### 3.2.1. Displaced Block-Tridiagonal Elimination Methods

Decomposition of $\Gamma^{D}$ leads to two other block tridiagonal elimination methods, that cannot be generated from $\Gamma^{A}$.

### Displaced Block-Tridiagonal Row (DBTR) Elimination■

$$\Gamma^{D}=\begin{bmatrix} & Ip & \\ \hdashline Bn^{e} & An^{\circ} \\ Bm^{e} & Am^{\circ} \end{bmatrix}\left[\begin{array}{c:c} \breve{\sigma}^{D} & \begin{matrix} Dn \\ Dm \end{matrix} \end{array}\right]$$

The method introduces the two extraneous blocks $Bn^{e}$ and $Bm^{e}$. An efficient sequence of its calculations decomposes the stem $\sigma^{D}$ through block row manipulation as in BRBR elimination.

Stem: $\breve{B}n^{\#}$(A5), $Bm''$(A17), $Bm^{*}$(A18), $Cm^{\#}$(A3a), $\breve{C}m^{\#}$(A6)
Fins: $An'$(A9), $Am'$(A7), $An^{\circ}$(A10), $Am^{\circ}$(A8)
Head: $Dm^{\#}$(A4a), $Bn^{\#}$(A2c), $Bm^{\#}$(A1c)

Note that $Bm^{e}$ and $Bn^{e}$ are not evaluated, since they need not to be stored. Their contribution to the head is estimated through substituting their expressions of (A23) and (A24) into (A27) and (A28) leading to (A2c) and (A1c), respectively.

### Displaced Block-Tridiagonal Column (DBTC) Elimination

$$\Gamma^{D}=\begin{bmatrix} & \breve{\sigma}^{D} \\ \hdashline 0n & An \\ & Am \end{bmatrix}\left[\begin{array}{c:c} Ip & \begin{matrix} Dn^{*} \\ Dm^{*} \end{matrix} \end{array}\right]$$

The method, as described in Ref. [11], decomposes the stem $\sigma^{D}$ through scalar column manipulation as in SCSC elimination

However, as with ABTR elimination, it would prove profitable to decompose the stem through an alternating manipulation method as in !C?R elimination, where ! and ? stand for S or B. Using SCSR elimination, the following sequence of manipulations applies.

Stem: $Ln\hat{U}n$(A5), $Cn'$(A11), $Bm''$(A17), $Cm^{\#}$(A3b), $\hat{L}mUm$(A6)
Fins: $Dn'$(A13), $Dm^{\#}$(A4b), $Dm''$(A19), $Dm^{*}$(A20), $Dn''$(A21), $Dn^{*}$(A22)
Head: $Bn^{\#}$(A2a), $Bm^{\#}$(A1a)

### 3.3 Solution Procedure

We observe that a decomposition $\mathbb{G}^{+}=\mathbb{L}\mathbb{U}^{+}$ that is based on any of the elimination methods described in §3.1 and §3.2 can be constructed from the corresponding decomposition $\Gamma_{j}^{+}=\mu_{j}\nu_{j}^{+}$ of the module. The significant part of $\mathbb{L}$ (or $\mathbb{U}^{+}$) is obtained by arranging the $\mu_{j}$'s (or $\nu_{j}^{+}$'s), sequentially, in order.

The procedure for solving the matrix equation $\mathbb{G}^{+}\mathsf{z}=\mathsf{g}$, which can be described in terms of manipulation of the augmented matrix $\mathbb{G}^{+}\equiv[\mathbb{G} \,\vdots\, \mathsf{g}]$, can, similarly, be described in terms of manipulation of the augmented module $\Gamma_j^+ \equiv [\Gamma_j \,\vdots\, \gamma_j]$. The manipulation of $\mathbb{G}^+$ applies a forward sweep which corresponds to the decomposition $\mathbb{G}^+ = \mathbb{L}\mathbb{U}^+ \equiv \mathbb{L}[\mathbb{U} \,\vdots\, Z]$, that is followed by a backward sweep which corresponds to the decomposition $\mathbb{U}^+ = \mathbb{U}[\mathrm{I} \,\vdots\, \mathrm{z}]$. Similarly, the manipulation of $\Gamma_j^+$ applies a forward sweep involving two steps. The first step performs the decomposition $\Gamma_j^{\;+} = \mu_j \nu_j^+$. The second step evaluates the head $\eta_j^{\;+} = \Psi_j \Phi_j^+$, then applies the head-tail relation to determine the tail of $\Gamma_{j+1}^+$. In a backward sweep, two steps are applied to $\nu_j^+ \equiv [\mathrm{N}_j \quad \Phi_j \,\vdots\, \mathrm{Z}_j]$ leading to the solution module $z_j$ ($z_j^A = [zm_j^{\;t} \quad zn_j^{\;t}]^t$, $z_j^D = [zn_j^{\;t} \quad zm_{j+1}^{\;t}]^t$). With $z_{j+1}$ known, the first step uses $\xi_{j+1}$ ($\xi_{j+1}^A = z_{j+1}^A$, $\xi_{j+1}^D = zn_{j+1}$) in the back substitution relation $\mathrm{Z}_j - \Phi_j \xi_{j+1} = \zeta_j$ to contract $\nu_j^+$ to $\mathrm{N}_j^+ \equiv [\mathrm{N}_j \,\vdots\, \zeta_j]$. The second step solves $\mathrm{N}_j z_j = \zeta_j$ for $z_j$ which is equivalent to the decomposition $\mathrm{N}_j^+ = \mathrm{N}_j [Ip \,\vdots\, z_j]$.

## 4. OPERATION COUNTS AND STORAGE NEEDS

The operation counts are measured by the number of multiplications (mul) with the understanding that the number of additions is comparable. The storage needs are measured by the number of locations (loc) required to store arrays calculated in the forward sweep for use in the backward sweep; provided that the elements of $\mathbb{G}^+$ are not stored but are generated when needed, as is usually done in the numerical solution of a set of ODEs, for example.

The modules introduced in §3 allow easy evaluation of the elimination methods. Per module (i.e., per grid point), each method requires, for the manipulation of $\mathbb{G}^+$, as many operations as it requires to manipulate $\Gamma^+$, and as many storage locations as it requires to store $\nu^+ \equiv [\nu \,\vdots\, \mathrm{Z}]$.

The manipulation of $\Gamma_j^+ \equiv [\Gamma_j \,\vdots\, \gamma_j]$ involves manipulating $\Gamma_j$ which requires the decomposition of the stem $\sigma_j$, the processing of the fins $\psi_j$ and $\phi_j$, and the handling of the head $\eta_j$; as well as manipulating $\gamma_j$ which requires the evaluation of $\mathrm{Z}_j$ then $\theta_j$ in the forward sweep, and $\zeta_j$ then $z_j$ in the backward sweep.

All methods require $\dfrac{p^3 - p}{3}$ (mul) for decomposing the stem. To illustrate, a scalar decomposition of the stem $\sigma_j^A$ involves the use of Eqs. (A6,19,9,2,5) which requires $\dfrac{m^3 - m}{3}, \dfrac{m^2 \pm m}{2} n, \dfrac{m^2 \mp m}{2} n, mn^2, \dfrac{n^3 - n}{3}$ (mul), respectively, totaling $\dfrac{p^3 - p}{3}$ (mul); whereas a scalar decomposition of the stem $\sigma_j^D$ involves the use of Eqs. (A5,11,17,3,6) which requires $\dfrac{n^3 - n}{3}, \dfrac{n^2 \mp n}{2} m, \dfrac{n^2 \pm n}{2} m, m^2 n, \dfrac{m^3 - m}{3}$ (mul), respectively, totaling, again, $\dfrac{p^3 - p}{3}$ (mul). Scalar/Scalar elimination methods apply such scalar decomposition of $\sigma_j$.

BRBR\BCBC\BCBR\BRBC elimination applies block decomposition of $\sigma_j^A$ or $\sigma_j^D$ which involves the use of Eqs. (A6,9\19\19\9,10\20\20\10, 2,5) or (A5,17\11\17\11, 18\12\18\12,3,6), requiring, respectively, $\frac{m^3-m}{3}, \frac{m^2 \pm m}{2}n, \frac{m^2 \mp m}{2}n, mn^2, \frac{n^3-n}{3}$ (mul) or $\frac{n^3-n}{3}, \frac{n^2 \mp n}{2}m, \frac{n^2 \pm n}{2}m, m^2n, \frac{m^3-m}{3}$ (mul), totaling in either case $\frac{p^3-p}{3}$ (mul). ABT\DBT row or column elimination applies an identity decomposition involving either scalar or block decomposition of $\sigma_j^A \backslash \sigma_j^D$.

All methods, which introduce no extraneous elements, require *pmn* (mul) for evaluating the head $\eta_j = \Psi_j \Phi_j$ in preparation for evaluating the tail $\tau_{j+1}$ with the help of the head-tail relation, which involves no multiplications. Moreover, all such methods require $2p^2r$ (mul) to handle the right module $\gamma_j$: The determination of the head $\theta_j = \Psi_j \mathrm{Z}_j$ in the forward sweep requires *pmr* (mul), while the back substitution relation $\mathrm{Z}_j - \Phi_j \xi_{j+1} = \zeta_j$ in the backward sweep requires *pnr* (mul); adding up to $p^2r$ (mul). The solution of $\mathrm{M}_j \mathrm{Z}_j = \rho_j$ for $\mathrm{Z}_j$ in the forward sweep, and the solution of $\mathrm{N}_j z_j = \zeta_j$ for $z_j$ in the backward sweep, together require $p^2r$ (mul), irrespective of the decomposition $\mathrm{M}_j \mathrm{N}_j = \breve{\sigma}_j$ of the stem. These methods differ only in the operation counts for processing the fins $\psi_j$ and $\phi_j$, with BCR elimination requiring the least count *pmn* (mul).

The only two methods, ABTC and DBTR elimination, which introduce extraneous nonzero elements through elimination, perform additional calculations for evaluating the extraneous blocks and/or the blocks they contribute to.

As for the storage needs, all methods require *pr* (loc) to store $\mathrm{Z}_j$. They differ in the number of locations needed for storing the $\nu_j$'s presented in §3, with DBTC elimination requiring the least number *pn* (loc). Note that, in methods involving scalar elimination, square blocks need to be reserved for storing the triangular blocks *Um* and/or *Un*.

Table (1B) of Appendix B contains the above information, allowing for clear comparison among the methods. It is noted that some methods show preference to *m*>*n* and others to *m*<*n*. Taking *m*>*n*, a fair comparison is possible by interchanging *m* and *n* for the methods which prefer *m*<*n*. (When solving a set of ODEs, this amounts to assigning *j*=1 to the side with *n* boundary conditions.) The four methods SCSR, SRSC, BCBR, and BRBC whose operation counts are not biased toward *m* or *n* can be applied on two processors of a parallel machine; one operating from *j*=1 and the other from *j*=*J*, with equal efficiency, considerably reducing the calculation time.

The savings achieved by BCBR elimination in operations and by DBTC elimination in storage are of leading order significance when *p*>>1, in the two distinguished limits *m*~*n*~*p*/2 and (*m*,*n*)~(*p*,1). In comparison with the scalar elimination methods, savings of at least $\sim p^3/8$ (mul) and $\sim p^2/4$ (loc) are achieved in the former limit, while savings of at least $\sim p^3/2$ (mul) and $\sim p^2/2$ (loc) are achieved in the latter limit.

## 5. PIVOTING

The role of pivoting is to guarantee validity and stability. Partial pivoting (aiming at evading division by a zero) is all one needs to insure validity. It is also generally accepted that partial pivoting (aiming further at binding the multipliers) leads to stability [12]. Partial pivoting strategies involve either one-dimensional search for the pivot (e.g., row pivoting, column pivoting, and Lam's alternating pivoting) or the more expensive two-dimensional search (e.g., Keller's mixed pivoting [10] and the to-be-introduced local pivoting).

Lam's alternating pivoting [8] is given, here, special attention. It involves one-dimensional search for the pivot, but unlike uni-directional (row or column) pivoting, it introduces no extraneous elements and, furthermore, it guarantees validity. Expressed in terms of the present notation, Lam's strategy applies column pivoting to $\tau^A = \left[ Cm^\# \quad Dm^\# \right]$ in order to form and decompose a nonsingular pivotal block $Cm^\#$, and applies row pivoting to $\tau^D = \left[ Bn^{\# t} \quad Bm^{\# t} \right]^t$ in order to form and decompose a nonsingular pivotal block $Bn^\#$. These are valid processes since $\tau^A$ is of rank *m* and $\tau^D$ is of rank *n*, as can be shown following the reasoning of Keller [10, §5, Theorem]. It is emphasized here that applying column interchanges within $\tau^A$ is mandatory for guaranteed validity. Row interchanges among the *m* linearly independent rows of $\tau^A$ do not guarantee validity as the first *m* columns may be linearly dependent. Likewise, row interchanges within $\tau^D$ are mandatory.

As advocated by Varah [9], Lam's alternating pivoting should be followed by correspondingly alternating elimination. Step by step, column pivoting is to be followed by column elimination and row pivoting by row elimination. This enhances stability since the multipliers are then bounded by unity. Obviously, SCSR, SCBR, BCSR and BCBR elimination methods can benefit from this. So can do ABTR and DBTC elimination methods- contrary to the common belief- if the decomposition of the stem $\sigma$ is carried out as in SCSR, SCBR, BCSR or BCBR elimination.

To enhance stability further we introduce the *Local Pivoting Strategy* which applies full pivoting (maximum pivot strategy) to the segments $\tau^A$ and $\tau^D$. Only those methods which inflect *Dm* to at least $Dm''$ and *Bm* to at least $Bm''$ can apply local pivoting. (Note that Keller's mixed pivoting [10], if interpreted as applying full pivoting to $\tau^A$ and row pivoting to $\tau^D$, is midway between Lam's and local pivoting.)

Table (1B) indicates the admissibility of Lam's and/or local pivoting by an asterisk *.

## 6. BORDERED AND DOUBLE-BORDERED SYSTEMS

Extension of modular analysis to bordered (BABD) and double-bordered (DBABD) almost block diagonal systems of equations is possible. For example, for a DBABD system with right-bordered array of *q* columns

$$\left[ Hm_1{}^t \quad Hn_1{}^t \quad \cdots \quad Hm_j{}^t \quad Hn_j{}^t \quad \cdots \quad Hm_J{}^t \quad Hn_J{}^t \quad Hk_1{}^t \right]^t$$

and lower-bordered array of *k* ($\leq n+q$) rows

$$\left[ Ak_1 \quad Bk_1 \quad \cdots \quad Ak_j \quad Bk_j \quad \cdots \quad Ak_J \quad Bk_J \quad Hk_1 \quad gk_1 \right]$$

we define the A-module

$$\Gamma_j^{+A} \equiv \left[\begin{array}{c:c} \sigma_j^A & \phi_j^{+A} \\ \hdashline \psi_j^A & \eta_j^{+A} \end{array}\right] = \left[\begin{array}{cc:cc:cc} Cm_j^{\#} & Dm_j^{\#} & & & Hm_j^{\#} & gm_j^{\#} \\ An_j & Bn_j & Cn_{j+1} & Dn_{j+1} & Hn_j & gn_j \\ \hdashline Am_j & Bm_j & Cm_{j+1}^{\Rightarrow} & Dm_{j+1}^{\Rightarrow} & Hm_{j+1}^{\Rightarrow} & gm_{j+1}^{\Rightarrow} \\ \hdashline Ak_j^{\#} & Bk_j^{\#} & Ak_{j+1}^{\Rightarrow} & Bk_{j+1}^{\Rightarrow} & Hk_{j+1}^{\Rightarrow} & gk_{j+1}^{\Rightarrow} \end{array}\right]$$

Its tail

$$\tau_j^{+A} = \left[\begin{array}{cc:cc} Cm_j^{\#} & Dm_j^{\#} & Hm_j^{\#} & gm_j^{\#} \\ \hdashline Ak_j^{\#} & Bk_j^{\#} & Hk_j^{\#} & gk_j^{\#} \end{array}\right]$$

having the hidden segment $\left[Hk_j^{\#} \quad gk_j^{\#}\right]$ which does not appear in $\Gamma_j^{+A}$, is determined through the head-tail relation

$$\eta_{j-1}^{+A} + \tau_j^{+A} = \left[\begin{array}{cc:cc} Cm_j & Dm_j & Hm_j & gm_j \\ \hdashline Ak_j & Bk_j & Hk_{j-1}^{\#} & gk_{j-1}^{\#} \end{array}\right]$$

All twenty elimination methods given in §3 are applicable. Handling the new bordered arrays requires $pq$ (loc) but different numbers of (mul). The methods that do not introduce extraneous elements require $2p^2q + pq^2 + 2pqr$ (mul), when $k=q$.

## 6. CONCLUSION

Using the novel approach of modular analysis, we have analyzed the sequential solution methods for almost block diagonal systems of equations. Two modules have been identified and have made it possible to express and assess all possible band and block elimination methods on the basis of their operation counts, storage needs, and admissibility of partial pivoting. Modular analysis has also been extended to bordered almost block diagonal systems.

Implementation of the four distinguished methods SCSR, BCSR, BCBR and DBTC within the COLROW platform [13], which was designed to give SCSR elimination its best performance, showed that the other three methods could outperform SCSR, in some cases (when *m*>>*n*). FORTRAN codes and related materials are given in Appendix C.

**Appendix A: Inflections of Blocks of $\mathbb{G}$**

In §3, two modules, $\Gamma^A$ and $\Gamma^S$, of the matrix of coefficients $\mathbb{G}$ are introduced and decomposed to generate the elimination methods. The process involves inflections of the blocks of $\mathbb{G}$, which proceed for a block *E*, say, according to the following scheme.

$$\begin{array}{ccccccc} & & E & \rightarrow & E' & \rightarrow & E^\circ \\ & & \downarrow & & \downarrow & & \\ \breve{E}^\# & \xleftarrow{\text{if pivota}} & E^\# & \rightarrow & E'' & \rightarrow & E* \end{array}$$

The following equalities are to be used to determine the blocks with underscored leading character. The number of multiplications involved is given between wiggly brackets. When a $\pm$ or $\mp$ sign appears, the upper (or lower) sign corresponds to the triangular matrix *U* (or *L*) involved being unit diagonal.

$E^{\#}$ - Blocks

$$Bm - \underline{B}m^\# \overset{a}{=} AmDm* \overset{b}{=} Am'Dm'' \overset{c}{=} Am^\circ Dm^\# \quad \{\, m^2 n \,\} \tag{A1}$$

$$Bn - \underline{B}n^\# \overset{a}{=} AnDm* \overset{b}{=} An'Dm'' \overset{c}{=} An^\circ Dm^\# \quad \{\, mn^2 \,\} \tag{A2}$$

$$Cm - \underline{C}m^\# \overset{a}{=} Bm*Cn \overset{b}{=} Bm''Cn' \overset{c}{=} Bm^\# Cn^\circ \quad \{\, m^2 n \,\} \tag{A3}$$

$$Dm - \underline{D}m^\# \overset{a}{=} Bm*Dn \overset{b}{=} Bm''Dn' \overset{c}{=} Bm^\# Dn^\circ \quad \{\, mn^2 \,\} \tag{A4}$$

$\breve{E}^{\#}$ - Blocks (Decomposed blocks)

$$Bn^{\#} = \underline{\breve{B}}n^{\#} (\equiv \underline{L}n\underline{U}n) \ \{\frac{n^3-n}{3}\} \quad Cm^{\#} = \underline{\breve{C}}m^{\#} (\equiv \underline{L}m\underline{U}m) \ \{\frac{m^3-m}{3}\} \qquad \text{(A5,6)}$$

$E'$ and $E^{\circ}$ - Blocks

$$\underline{A}m'Um = Am \quad \{\frac{m^2 \mp m}{2} m\} \qquad \underline{A}m^{\circ}Lm = Am' \quad \{\frac{m^2 \pm m}{2} m\} \qquad \text{(A7,8)}$$

$$\underline{A}n'Um = An \quad \{\frac{m^2 \mp m}{2} n\} \qquad \underline{A}n^{\circ}Lm = An' \quad \{\frac{m^2 \pm m}{2} n\} \qquad \text{(A9,10)}$$

$$Ln\underline{C}n' = Cn \quad \{\frac{n^2 \pm n}{2} m\} \qquad Un\underline{C}n^{\circ} = Cn' \quad \{\frac{n^2 \mp n}{2} m\} \qquad \text{(A11,12)}$$

$$Ln\underline{D}n' = Dn \quad \{\frac{n^2 \pm n}{2} n\} \qquad Un\underline{D}n^{\circ} = Dn' \quad \{\frac{n^2 \mp n}{2} n\} \qquad \text{(A13,14)}$$

$E''$ and $E*$-Blocks

$$\underline{A}m'' = Am' - Bm*An' \quad \{m^2 n\} \quad \underline{A}m*Lm = Am'' \quad \{\frac{m^2 \pm m}{2} m\} \qquad \text{(A15,16)}$$

$$\underline{B}m''Un = Bm^{\#} \quad \{\frac{n^2 \mp n}{2} m\} \qquad \underline{B}m*Ln = Bm'' \quad \{\frac{n^2 \pm n}{2} m\} \qquad \text{(A17,18)}$$

$$Lm\underline{D}m'' = Dm^{\#} \quad \{\frac{m^2 \pm m}{2} n\} \qquad Um\underline{D}m* = Dm'' \quad \{\frac{m^2 \mp m}{2} n\} \qquad \text{(A19,20)}$$

$$\underline{D}n'' = Dn' - Cn'Dm* \quad \{mn^2\} \qquad Un\underline{D}n* = Dn'' \quad \{\frac{n^2 \mp n}{2} n\} \qquad \text{(A21,22)}$$

$E^{e}$ - Blocks and related $E^{\#}$ - Blocks

$$\underline{B}m^{e} = -Am^{\circ}Bm* \quad \{m^2 n\} \qquad \underline{B}n^{e} = -An^{\circ}Bm* \quad \{mn^2\} \qquad \text{(A23,24)}$$

$$\underline{C}m^{e} = -Dm*Cn^{\circ} \quad \{m^2 n\} \qquad Um\underline{C}m^{e} = -Dm''Cn^{\circ} \quad \{m^2 n + \frac{m^2 \mp m}{2} m\} \qquad \text{(A25a,b)}$$

$$\underline{D}m^{e} = -Dm*Dn^{\circ} \quad \{mn^2\} \qquad Um\underline{D}m^{e} = -Dm''Dn^{\circ} \quad \{mn^2 + \frac{m^2 \mp m}{2} n\} \qquad \text{(A26a,b)}$$

$$Bm - \underline{B}m^{\#} = Bm^{e}Dn + Am^{\circ}Dm \quad \{mn^2 + m^2 n\} \qquad \text{(A27)}$$

$$Bn - \underline{B}n^{\#} = Bn^{e}Dn + An^{\circ}Dm \quad \{n^3 + mn^2\} \qquad \text{(A28)}$$

$$Cm - \underline{C}m^{\#} = AmCm^{e} + BmCn^{\circ} \quad \{m^3 + m^2 n\} \qquad \text{(A29)}$$

$$Dm - \underline{D}m^{\#} = AmDm^{e} + BmDn^{\circ} \quad \{m^2 n + mn^2\} \qquad \text{(A30)}$$

## Appendix B: Assessment Table

The operation counts, storage needs, and admissibility of Lam's and local pivoting of all twenty elimination methods are contained in Table (B1) given in this appendix.

Note that the methods involving scalar elimination reserve square blocks for storing the triangular blocks *Um* and/or *Un*. The expressions given in the footnotes of Table (B1) should, therefore, replace the corresponding table entries.

Table (B1) : Assessment of Elimination Methods

| Method | Operation Counts (mul) | Storage Needs (loc) | Pivoting | |
|---|---|---|---|---|
| | $\frac{1}{3}\left(p^3 - p\right) + 2pmn + 2p^2r +$ | $pr + pn \; +$ | Lam’s | Local |
| SRSR | $\frac{1}{2}\left(m^3 + n^3 + m^2 - n^2\right)$ | $\frac{1}{2}\left(p^2 + m + n\right)^{a}$ | | * |
| SCSC | $\frac{1}{2}\left(m^3 + n^3 - m^2 + n^2\right)$ | $\frac{1}{2}\left(p^2 - m - n\right)^{a}$ | | * |
| SCSR | $\frac{1}{2}\left(m^3 + n^3 - m^2 - n^2\right)$ | $\frac{1}{2}\left(p^2 - m + n\right)^{a}$ | * | * |
| SRSC■ | $\frac{1}{2}\left(m^3 + n^3 + m^2 + n^2\right)$ | $\frac{1}{2}\left(p^2 + m - n\right)^{a}$ | | * |
| SRBR■ | $\frac{1}{2}\left(m^3 + m^2\right)$ | $\frac{1}{2}\left(p^2 + n^2 + m\right)^{a}$ | | * |
| SCBR■ | $\frac{1}{2}\left(m^3 - m^2\right)$ | $\frac{1}{2}\left(p^2 + n^2 - m\right)^{a}$ | * | * |
| SRBC■ | $\frac{1}{2}\left(m^3 + 2n^3 + m^2\right)$ | $\frac{1}{2}\left(p^2 - n^2 + m\right)^{b}$ | | |
| SCBC■ | $\frac{1}{2}\left(m^3 + 2n^3 - m^2\right)$ | $\frac{1}{2}\left(p^2 - n^2 - m\right)^{b}$ | | |
| BRSR■ | $\frac{1}{2}\left(2m^3 + n^3 + n^2\right)$ | $\frac{1}{2}\left(p^2 + m^2 + n\right)^{a}$ | | |
| BRSC■ | $\frac{1}{2}\left(2m^3 + n^3 - n^2\right)$ | $\frac{1}{2}\left(p^2 + m^2 - n\right)^{a}$ | | |
| BCSR■ | $\frac{1}{2}\left(n^3 - n^2\right)$ | $\frac{1}{2}\left(p^2 - m^2 + n\right)^{c}$ | * | * |
| BCSC■ | $\frac{1}{2}\left(n^3 + n^2\right)$ | $\frac{1}{2}\left(p^2 - m^2 - n\right)^{c}$ | | * |
| BRBR■ | $m^3$ | $p^2 - mn$ | | |
| BCBC■ | $n^3$ | $mn$ | | |
| BCBR■ | 0 | $pn$ | * | * |
| BRBC■ | $m^3 + n^3$ | $pm$ | | |
| ABTR | $pm^2$ | $p^2$ | * | * |
| DBTC | $pn^2$ | 0 | * | * |
| DBTR■ | $m^3 + (p+n)mr$ | $p^2$ | | |
| ABTCe | $pmn + n^3 + pmr$ | $pm$ | | |
| ABTCt | $p\left(\frac{3}{2}m^2 + n^2 + \frac{1}{2}m\right) + pmr$ | $pm$ | | |

[a] $p^2 - mn$

[b] $pm$

[c] $pn$

## Appendix C: Supplementary Materials

The following supplementary materials include listings of four FORTRAN codes:
SCSR.FOR; a Fortran program for the SCSR elimination method,
BCSR.FOR; a Fortran program for the BCSR elimination method,
BCBR.FOR; a Fortran program for the BCBR elimination method,
DBTC.FOR; a Fortran program for the DBTC elimination method.

SCSR.FOR is essentially the COLROW algorithm, without FLAG statements and with renaming of some variables.

BCSR.FOR, BCBR.FOR, and DBTC.FOR are modifications of SCSR.FOR which apply BCSR, BCBR, and DBTC eliminations, instead of SCSR elimination.

In the four codes, the integer variables:
LM is the number of rows in the top block,
LN is the number of rows in the bottom block,
KK is the number of repetitions of the repeated block, and
MM is the number of times the solution is to be carried out.
They appear in PARAMETER statements in the main program and in the two subroutines, where they currently assume the values LM=10, LN=1, KK=10, MM=1000000. A change in the values of LM and LN should be applied to all PARAMETER statements, keeping LM+LN=11.

The supplementary materials also include two input data listings:
TESTA; containing the top and bottom blocks, and
TESTB; containing the repeated block.

The data in TESTA and TESTB correspond to a coefficient matrix whose entries are chosen randomly. The right hand side column is calculated such that the solution column is $(1.0\ \ 1.1\ \ 1.2\ \ \ldots\ 1.8\ \ 1.9\ \ 2.0)^t$ repeated KK+1 times; irrespective of the values of LM and LN. A change in the values of LM and LN can be effected by positioning the line
BOTTOM BLOCK
in TESTA; according to the new values. Currently, it is positioned such that LM=10, LN=1.

## SCSR.FOR Fortran program

```
CCCCC
C    SCSR
CCCCC
      IMPLICIT REAL*8(A-H,O-Z)
      PARAMETER(LM=10,LN=1,KK=10,MM=1000000,
     ,  LM1=LM+1,LN1=LN+1,KK1=KK+1,
     ,  LC=LM+LN,LC1=LC+1,LA=2*LC,LA1=LA+1,N=LC*KK1)
      DIMENSION TOPBLK(LM,LC),ARRAY(LC,LA,KK),BOTBLK(LN,LC)
      DIMENSION IPIVOT(N),B(N)
      DIMENSION A1(LM,LC1),A2(LN,LC1),A3(LC,LA1)
      OPEN(1,FILE='TESTA')
      READ(1,1000) ((A1(I,J),J=1,LC1),I=1,LM)
      READ(1,1000) ((A2(I,J),J=1,LC1),I=1,LN)
      CLOSE(1)
      OPEN(2,FILE='TESTB')
      READ(2,2000) ((A3(I,J),J=1,LA1),I=1,LC)
      CLOSE(2)
      DO 4 M=1,MM
      DO 1 I=1,LM
      B(I)=A1(I,LC1)
      DO 1 J=1,LC
```

```
    1 TOPBLK(I,J)=A1(I,J)
      DO 2 K=1,KK
      DO 2 I=1,LC
      B(I+K*LC-LN)=A3(I,LA1)
      DO 2 J=1,LA
    2 ARRAY(I,J,K)=A3(I,J)
      DO 3 I=1,LN
      B(I+N-LN)=A2(I,LC1)
      DO 3 J=1,LC
    3 BOTBLK(I,J)=A2(I,J)
      CALL CRDCMP(TOPBLK,ARRAY,BOTBLK,IPIVOT,B)
      CALL CRSLVE(TOPBLK,ARRAY,BOTBLK,IPIVOT,B)
    4 CONTINUE
      OPEN(3,FILE='TESTC')
      WRITE (3,3000) (B(I),I=1,N)
      CLOSE(3)
      STOP
 1000 FORMAT(5X/(11F6.0,F8.0))
 2000 FORMAT(5X/(22F6.0,F8.0))
 3000 FORMAT(3E25.15)
      END
      SUBROUTINE CRDCMP(TOPBLK,ARRAY,BOTBLK,IPIVOT,B)
      IMPLICIT REAL*8(A-H,O-Z)
      PARAMETER(LM=10,LN=1,KK=10,
     ,  LM1=LM+1,LN1=LN+1,KK1=KK+1,
     ,  LC=LM+LN,LC1=LC+1,LA=2*LC,N=LC*KK1)
      DIMENSION TOPBLK(LM,LC),ARRAY(LC,LA,KK),BOTBLK(LN,LC)
      DIMENSION IPIVOT(N),B(N)
      LB=LC+LM
            DO 110 I=1,LM
        I1 = I + 1
        IPVT = I
        COLMAX = DABS(TOPBLK(I,I))
        DO 30 J=I1,LC
          TEMPIV = DABS(TOPBLK(I,J))
          IF (TEMPIV.LE.COLMAX) GO TO 30
          IPVT = J
          COLMAX = TEMPIV
   30   CONTINUE
        IPIVOT(I) = IPVT
        IF (IPVT.EQ.I) GO TO 60
        DO 40 L=I,LM
          SWAP = TOPBLK(L,IPVT)
          TOPBLK(L,IPVT) = TOPBLK(L,I)
          TOPBLK(L,I) = SWAP
   40   CONTINUE
        DO 50 L=1,LC
          SWAP = ARRAY(L,IPVT,1)
          ARRAY(L,IPVT,1) = ARRAY(L,I,1)
          ARRAY(L,I,1) = SWAP
   50   CONTINUE
   60   CONTINUE
        CII = TOPBLK(I,I)
        DO 80 J=I1,LC
          CIJ = TOPBLK(I,J)/CII
          TOPBLK(I,J) = CIJ
          DO 70 L=I1,LM
   70       TOPBLK(L,J) = TOPBLK(L,J) - TOPBLK(L,I)*CIJ
          DO 75 L=1,LC
   75       ARRAY(L,J,1) = ARRAY(L,J,1) - ARRAY(L,I,1)*CIJ
   80   CONTINUE
      BI=B(I)/CII
      B(I)=BI
      DO 90 L=I1,LM
   90 B(L)=B(L)-TOPBLK(L,I)*BI
      DO 95 L=LM1,LB
   95 B(L)=B(L)-ARRAY(L-LM,I,1)*BI
  110 CONTINUE
      INCR = 0
      DO 320 K=1,KK
        K1 = K + 1
          INCRM=INCR+LM
          INCRC=INCR+LC
          INCRB=INCR+LB
        DO 180 J=LM1,LC
```

```
        J1 = J + 1
        JN = J - LM
        JN1 = JN + 1
        INCRJ=INCR+J
        IPVT = JN
        ROWMAX = DABS(ARRAY(JN,J,K))
       DO 120 I=JN1,LC
          TEMPIV = DABS(ARRAY(I,J,K))
          IF (TEMPIV.LE.ROWMAX) GO TO 120
          IPVT = I
          ROWMAX = TEMPIV
120     CONTINUE
       INCRP=INCRM+IPVT
        IPIVOT(INCRJ)=INCRP
        IF (IPVT.EQ.JN) GO TO 140
        DO 130 L=J,LA
          SWAP = ARRAY(IPVT,L,K)
          ARRAY(IPVT,L,K) = ARRAY(JN,L,K)
          ARRAY(JN,L,K) = SWAP
130     CONTINUE
          SWAP = B(INCRP)
          B(INCRP) = B(INCRJ)
          B(INCRJ) = SWAP
140     CONTINUE
        ROWPIV = ARRAY(JN,J,K)
        DO 150 I=JN1,LC
          ARRAY(I,J,K) = ARRAY(I,J,K)/ROWPIV
150     CONTINUE
        DO 160 L=J1,LA
          ROWMLT = ARRAY(JN,L,K)
          DO 160 I=JN1,LC
160          ARRAY(I,L,K) = ARRAY(I,L,K) - ROWMLT*ARRAY(I,J,K)
          ROWMLT = B(INCRJ)
          DO 170 I=JN1,LC
            INCRI=INCRM+I
170          B(INCRI) = B(INCRI) - ROWMLT*ARRAY(I,J,K)
180   CONTINUE
      DO 310 I=LN1,LC
        IM = I + LM
        IM1 = IM + 1
        IMC = IM - LC
        INCRN = INCR + IM
        IN=I-LN
        I1 = I + 1
        IPVT = IM
        COLMAX = DABS(ARRAY(I,IPVT,K))
        DO 190 J=IM1,LA
          TEMPIV = DABS(ARRAY(I,J,K))
          IF (TEMPIV.LE.COLMAX) GO TO 190
          IPVT = J
          COLMAX = TEMPIV
190     CONTINUE
        IPIVOT(INCRN) = INCR + IPVT
        IF (IPVT.EQ.IM) GO TO 240
        DO 200 L=I,LC
          SWAP = ARRAY(L,IPVT,K)
          ARRAY(L,IPVT,K) = ARRAY(L,IM,K)
          ARRAY(L,IM,K) = SWAP
200     CONTINUE
        IPC = IPVT - LC
        IF (K.EQ.KK) GO TO 220
        DO 210 L=1,LC
          SWAP = ARRAY(L,IPC,K1)
          ARRAY(L,IPC,K1) = ARRAY(L,IMC,K1)
          ARRAY(L,IMC,K1) = SWAP
210     CONTINUE
        GO TO 240
220     CONTINUE
        DO 230 L=1,LN
          SWAP = BOTBLK(L,IPC)
          BOTBLK(L,IPC) = BOTBLK(L,IMC)
          BOTBLK(L,IMC) = SWAP
230     CONTINUE
240     CONTINUE
        CII = ARRAY(I,IM,K)
```

```
         DO 270 J=IM1,LA
           CIJ = ARRAY(I,J,K)/CII
           ARRAY(I,J,K) = CIJ
           DO 250 L=I1,LC
250          ARRAY(L,J,K) = ARRAY(L,J,K) - ARRAY(L,IM,K)*CIJ
           JC = J - LC
           IF (K.LT.KK) THEN
           DO 260 L=1,LC
260          ARRAY(L,JC,K1) = ARRAY(L,JC,K1) -ARRAY(L,IMC,K1)*CIJ
           ELSE
           DO 265 L=1,LN
265          BOTBLK(L,JC) = BOTBLK(L,JC) -BOTBLK(L,IMC)*CIJ
           ENDIF
270        CONTINUE
     BI=B(INCRN)/CII
     B(INCRN)=BI
     DO 280 L=I1,LC
     LINCRM=L+INCRM
280 B(LINCRM)=B(LINCRM)-ARRAY(L,IM,K)*BI
     IF (K.LT.KK) THEN
     DO 290 L=1,LC
     LINCRB=L+INCRB
290 B(LINCRB)=B(LINCRB)-ARRAY(L,IN,K1)*BI
     ELSE
     DO 295 L=1,LN
     LINCRB=L+INCRB
295 B(LINCRB)=B(LINCRB)-BOTBLK(L,IN)*BI
     ENDIF
310   CONTINUE
       INCR = INCRC
320 CONTINUE
     IF (LN.EQ.1) GO TO 400
     INCRM=INCR+LM
     LC2=LC-1
     DO 390 J=LM1,LC2
       J1 = J + 1
       JN2 = J - LM
       JN21 = JN2 + 1
       INCRJ=INCR+J
       IPVT = JN2
       ROWMAX = DABS(BOTBLK(JN2,J))
       DO 330 I=JN21,LN
         TEMPIV = DABS(BOTBLK(I,J))
         IF (TEMPIV.LE.ROWMAX) GO TO 330
         IPVT = I
         ROWMAX = TEMPIV
330   CONTINUE
       INCRP=INCRM+IPVT
       IPIVOT(INCRJ)=INCRP
       IF (IPVT.EQ.JN2) GO TO 350
       DO 340 L=J,LC
         SWAP = BOTBLK(IPVT,L)
         BOTBLK(IPVT,L) = BOTBLK(JN2,L)
         BOTBLK(JN2,L) = SWAP
340   CONTINUE
           SWAP = B(INCRP)
           B(INCRP) = B(INCRJ)
           B(INCRJ) = SWAP
350   CONTINUE
       ROWPIV = BOTBLK(JN2,J)
       DO 360 I=JN21,LN
360      BOTBLK(I,J) = BOTBLK(I,J)/ROWPIV
       DO 370 L=J1,LC
         ROWMLT = BOTBLK(JN2,L)
         DO 370 I=JN21,LN
370        BOTBLK(I,L) = BOTBLK(I,L) - ROWMLT*BOTBLK(I,J)
           ROWMLT = B(INCRJ)
           DO 380 I=JN21,LN
             INCRI=INCRM+I
380          B(INCRI) = B(INCRI) - ROWMLT*BOTBLK(I,J)
390 CONTINUE
400 CONTINUE
     RETURN
     END
     SUBROUTINE CRSLVE(TOPBLK,ARRAY,BOTBLK,IPIVOT,B)
```

```
      IMPLICIT REAL*8(A-H,O-Z)
      PARAMETER(LM=10,LN=1,KK=10,
     ,  LM1=LM+1,LN1=LN+1,KK1=KK+1,
     ,  LC=LM+LN,LC1=LC+1,LA=2*LC,N=LC*KK1)
      DIMENSION TOPBLK(LM,LC),ARRAY(LC,LA,KK),BOTBLK(LN,LC)
      DIMENSION IPIVOT(N),B(N)
      INCR=KK*LC
      INCRM=INCR+LM
      DO 210 L=1,LN
        J = LC1 - L
        INCRJ = INCR + J
        LN1L = LN1 - L
        B(INCRJ) = B(INCRJ)/BOTBLK(LN1L,J)
        IF (L.EQ.LN) GO TO 200
        BINCRJ = B(INCRJ)
        LNL = LN - L
        DO 190 I=1,LNL
          INCRI = INCRM + I
          B(INCRI) = B(INCRI) - BOTBLK(I,J)*BINCRJ
190     CONTINUE
200     CONTINUE
210   CONTINUE
      DO 300 LK=1,KK
        K = KK1 - LK
        INCR = INCR - LC
        DO 240 L1=LN1,LC
          I = LC + LN1 - L1
          IM = I + LM
          IM1 = IM + 1
          INCRN = INCR + IM
          DOTPRD = B(INCRN)
          DO 220 J=IM1,LA
            INCRJ = INCR + J
            DOTPRD = DOTPRD - ARRAY(I,J,K)*B(INCRJ)
220       CONTINUE
          B(INCRN) = DOTPRD
          IPVTN = IPIVOT(INCRN)
          IF (INCRN.EQ.IPVTN) GO TO 230
          SWAP = B(INCRN)
          B(INCRN) = B(IPVTN)
          B(IPVTN) = SWAP
230       CONTINUE
240     CONTINUE
        INCRM = INCR + LM
        DO 260 J=LC1,LA
          INCRJ = INCR + J
          BINCRJ = B(INCRJ)
          DO 250 I=1,LN
            INCRI = INCRM + I
            B(INCRI) = B(INCRI) - ARRAY(I,J,K)*BINCRJ
250       CONTINUE
260     CONTINUE
        DO 290 L=1,LN
          J = LC1 - L
          INCRJ = INCR + J
          LN1L = LN1 - L
          B(INCRJ) = B(INCRJ)/ARRAY(LN1L,J,K)
          IF (L.EQ.LN) GO TO 280
          BINCRJ = B(INCRJ)
          LNL = LN - L
          DO 270 I=1,LNL
            INCRI = INCRM + I
            B(INCRI) = B(INCRI) - ARRAY(I,J,K)*BINCRJ
270       CONTINUE
280       CONTINUE
290     CONTINUE
300   CONTINUE
      DO 330 L=1,LM
        I = LM1 - L
        I1 = I + 1
        DOTPRD = B(I)
        DO 310 J=I1,LC
          DOTPRD = DOTPRD - TOPBLK(I,J)*B(J)
310     CONTINUE
        B(I) = DOTPRD
```

```
      IPVTI = IPIVOT(I)
      IF (I.EQ.IPVTI) GO TO 320
      SWAP = B(I)
      B(I) = B(IPVTI)
      B(IPVTI) = SWAP
 320   CONTINUE
 330 CONTINUE
     RETURN
     END
```

## BCSR.FOR Fortran program

```
CCCCC
C    BCSR
CCCCC
     IMPLICIT REAL*8(A-H,O-Z)
     PARAMETER(LM=10,LN=1,KK=10,MM=1000000,
    ,  LM1=LM+1,LN1=LN+1,KK1=KK+1,
    ,  LC=LM+LN,LC1=LC+1,LA=2*LC,LA1=LA+1,N=LC*KK1)
     DIMENSION TOPBLK(LM,LC),ARRAY(LC,LA,KK),BOTBLK(LN,LC)
     DIMENSION IPIVOT(N),B(N)
     DIMENSION A1(LM,LC1),A2(LN,LC1),A3(LC,LA1)
     OPEN(1,FILE='TESTA')
     READ(1,1000) ((A1(I,J),J=1,LC1),I=1,LM)
     READ(1,1000) ((A2(I,J),J=1,LC1),I=1,LN)
     CLOSE(1)
     OPEN(2,FILE='TESTB')
     READ(2,2000) ((A3(I,J),J=1,LA1),I=1,LC)
     CLOSE(2)
     DO 4 M=1,MM
     DO 1 I=1,LM
     B(I)=A1(I,LC1)
     DO 1 J=1,LC
   1 TOPBLK(I,J)=A1(I,J)
     DO 2 I=1,LN
     B(I+N-LN)=A2(I,LC1)
     DO 2 J=1,LC
   2 BOTBLK(I,J)=A2(I,J)
     DO 3 K=1,KK
     DO 3 I=1,LC
     B(I+K*LC-LN)=A3(I,LA1)
     DO 3 J=1,LA
   3 ARRAY(I,J,K)=A3(I,J)
     CALL CRDCMP(TOPBLK,ARRAY,BOTBLK,IPIVOT,B)
     CALL CRSLVE(TOPBLK,ARRAY,BOTBLK,IPIVOT,B)
   4 CONTINUE
     OPEN(3,FILE='TESTC')
     WRITE (3,3000) (B(I),I=1,N)
     CLOSE(3)
     STOP
 1000 FORMAT(5X/(11F6.0,F8.0))
 2000 FORMAT(5X/(22F6.0,F8.0))
 3000 FORMAT(3E25.15)
     END
     SUBROUTINE CRDCMP(TOPBLK,ARRAY,BOTBLK,IPIVOT,B)
     IMPLICIT REAL*8(A-H,O-Z)
     PARAMETER(LM=10,LN=1,KK=10,
    ,  LM1=LM+1,LN1=LN+1,KK1=KK+1,
    ,  LC=LM+LN,LC1=LC+1,LA=2*LC,N=LC*KK1)
     DIMENSION TOPBLK(LM,LC),ARRAY(LC,LA,KK),BOTBLK(LN,LC)
     DIMENSION IPIVOT(N),B(N)
     LB=LC+LM
     LB1=LB+1
     DO 80 I=1,LM
      I1 = I + 1
      IPVT = I
      COLMAX = DABS(TOPBLK(I,I))
      DO 30 J=I1,LC
        TEMPIV = DABS(TOPBLK(I,J))
        IF (TEMPIV.LE.COLMAX) GO TO 30
        IPVT = J
        COLMAX = TEMPIV
  30   CONTINUE
      IPIVOT(I) = IPVT
```

```
      IF (IPVT.EQ.I) GO TO 60
      DO 40 L=1,LM
        SWAP = TOPBLK(L,IPVT)
        TOPBLK(L,IPVT) = TOPBLK(L,I)
        TOPBLK(L,I) = SWAP
 40     CONTINUE
      DO 50 L=1,LC
        SWAP = ARRAY(L,IPVT,1)
        ARRAY(L,IPVT,1) = ARRAY(L,I,1)
        ARRAY(L,I,1) = SWAP
 50     CONTINUE
 60     CONTINUE
      CII = TOPBLK(I,I)
      DO 70 J=I1,LC
        CIJ = TOPBLK(I,J)/CII
        TOPBLK(I,J) = CIJ
        DO 70 L=I1,LM
 70       TOPBLK(L,J) = TOPBLK(L,J) - TOPBLK(L,I)*CIJ
    BI=B(I)/CII
    B(I)=BI
    DO 75 L=I1,LM
 75 B(L)=B(L)-TOPBLK(L,I)*BI
 80 CONTINUE
    DO 110 L=LM,1,-1
    L1=L-1
    BL=B(L)
    DO 90 I=1,L1
    CIL=TOPBLK(I,L)
    DO 85 J=LM1,LC
 85 TOPBLK(I,J)=TOPBLK(I,J)-CIL*TOPBLK(L,J)
 90 B(I)=B(I)-CIL*BL
    DO 100 I=1,LC
    IM=I+LM
    CIL=ARRAY(I,L,1)
    DO 95 J=LM1,LC
 95 ARRAY(I,J,1) = ARRAY(I,J,1) - CIL*TOPBLK(L,J)
100 B(IM)=B(IM)-CIL*BL
110 CONTINUE
    INCR = 0
    DO 320 K=1,KK
      K1 = K + 1
        INCRM=INCR+LM
        INCRC=INCR+LC
        INCRB=INCR+LB
      DO 180 J=LM1,LC
        J1 = J + 1
        JN = J - LM
        JN1 = JN + 1
        INCRJ=INCR+J
        IPVT = JN
        ROWMAX = DABS(ARRAY(JN,J,K))
        DO 120 I=JN1,LC
          TEMPIV = DABS(ARRAY(I,J,K))
          IF (TEMPIV.LE.ROWMAX) GO TO 120
          IPVT = I
          ROWMAX = TEMPIV
120       CONTINUE
        INCRP=INCRM+IPVT
        IPIVOT(INCRJ)=INCRP
        IF (IPVT.EQ.JN) GO TO 140
        DO 130 L=J,LA
          SWAP = ARRAY(IPVT,L,K)
          ARRAY(IPVT,L,K) = ARRAY(JN,L,K)
          ARRAY(JN,L,K) = SWAP
130       CONTINUE
          SWAP = B(INCRP)
          B(INCRP) = B(INCRJ)
          B(INCRJ) = SWAP
140       CONTINUE
        ROWPIV = ARRAY(JN,J,K)
        DO 150 I=JN1,LC
          ARRAY(I,J,K) = ARRAY(I,J,K)/ROWPIV
150       CONTINUE
        DO 160 L=J1,LA
          ROWMLT = ARRAY(JN,L,K)
```

```
         DO 160 I=JN1,LC
160          ARRAY(I,L,K) = ARRAY(I,L,K) - ROWMLT*ARRAY(I,J,K)
         ROWMLT = B(INCRJ)
         DO 170 I=JN1,LC
           INCRI=INCRM+I
170          B(INCRI) = B(INCRI) - ROWMLT*ARRAY(I,J,K)
180   CONTINUE
      DO 260 I=LN1,LC
        IM = I + LM
        IM1 = IM + 1
        IMC = IM - LC
        INCRN = INCR + IM
        IN=I-LN
        I1 = I + 1
        IPVT = IM
        COLMAX = DABS(ARRAY(I,IPVT,K))
        DO 190 J=IM1,LA
          TEMPIV = DABS(ARRAY(I,J,K))
          IF (TEMPIV.LE.COLMAX) GO TO 190
          IPVT = J
          COLMAX = TEMPIV
190     CONTINUE
        IPIVOT(INCRN) = INCR + IPVT
        IF (IPVT.EQ.IM) GO TO 240
        DO 200 L=LN1,LC
          SWAP = ARRAY(L,IPVT,K)
          ARRAY(L,IPVT,K) = ARRAY(L,IM,K)
          ARRAY(L,IM,K) = SWAP
200     CONTINUE
        IPC = IPVT - LC
        IF (K.EQ.KK) GO TO 220
        DO 210 L=1,LC
          SWAP = ARRAY(L,IPC,K1)
          ARRAY(L,IPC,K1) = ARRAY(L,IMC,K1)
          ARRAY(L,IMC,K1) = SWAP
210     CONTINUE
        GO TO 240
220     CONTINUE
        DO 230 L=1,LN
          SWAP = BOTBLK(L,IPC)
          BOTBLK(L,IPC) = BOTBLK(L,IMC)
          BOTBLK(L,IMC) = SWAP
230     CONTINUE
240     CONTINUE
        CII = ARRAY(I,IM,K)
        DO 250 J=IM1,LA
          CIJ = ARRAY(I,J,K)/CII
          ARRAY(I,J,K) = CIJ
          DO 250 L=I1,LC
250          ARRAY(L,J,K) = ARRAY(L,J,K) - ARRAY(L,IM,K)*CIJ
    BI=B(INCRN)/CII
    B(INCRN)=BI
    DO 255 L=I1,LC
    LINCRM=L+INCRM
255 B(LINCRM)=B(LINCRM)-ARRAY(L,IM,K)*BI
260 CONTINUE
    DO 300 L=LM,1,-1
    L1=L-1
    L1N=L1+LN
    BL=B(L+INCRC)
    DO 270 I=LN1,L1N
    IINCRM=I+INCRM
    CIL=ARRAY(I,L+LC,K)
    DO 265 J=LB1,LA
265 ARRAY(I,J,K)=ARRAY(I,J,K)-CIL*ARRAY(L+LN,J,K)
270 B(IINCRM)=B(IINCRM)-CIL*BL
    IF(K.LT.KK) THEN
    DO 280 I=1,LC
    IINCRB=I+INCRB
    CIL=ARRAY(I,L,K1)
    DO 275 J=LM1,LC
275 ARRAY(I,J,K1) = ARRAY(I,J,K1) - CIL*ARRAY(L+LN,J+LC,K)
280 B(IINCRB)=B(IINCRB)-CIL*BL
    ELSE
    DO 290 I=1,LN
```

```
   IINCRB=I+INCRB
   CIL=BOTBLK(I,L)
   DO 285 J=LM1,LC
285 BOTBLK(I,J) = BOTBLK(I,J) - CIL*ARRAY(L+LN,J+LC,K)
290 B(IINCRB)=B(IINCRB)-CIL*BL
   ENDIF
300 CONTINUE
     INCR = INCRC
320 CONTINUE
   IF (LN.EQ.1) GO TO 400
   INCRM=INCR+LM
   LC2=LC-1
   DO 390 J=LM1,LC2
     J1 = J + 1
     JN2 = J - LM
     JN21 = JN2 + 1
     INCRJ=INCR+J
     IPVT = JN2
     ROWMAX = DABS(BOTBLK(JN2,J))
     DO 330 I=JN21,LN
       TEMPIV = DABS(BOTBLK(I,J))
       IF (TEMPIV.LE.ROWMAX) GO TO 330
       IPVT = I
       ROWMAX = TEMPIV
330   CONTINUE
     INCRP=INCRM+IPVT
     IPIVOT(INCRJ)=INCRP
     IF (IPVT.EQ.JN2) GO TO 350
     DO 340 L=J,LC
       SWAP = BOTBLK(IPVT,L)
       BOTBLK(IPVT,L) = BOTBLK(JN2,L)
       BOTBLK(JN2,L) = SWAP
340   CONTINUE
         SWAP = B(INCRP)
         B(INCRP) = B(INCRJ)
         B(INCRJ) = SWAP
350   CONTINUE
     ROWPIV = BOTBLK(JN2,J)
     DO 360 I=JN21,LN
360     BOTBLK(I,J) = BOTBLK(I,J)/ROWPIV
     DO 370 L=J1,LC
       ROWMLT = BOTBLK(JN2,L)
       DO 370 I=JN21,LN
370       BOTBLK(I,L) = BOTBLK(I,L) - ROWMLT*BOTBLK(I,J)
         ROWMLT = B(INCRJ)
         DO 380 I=JN21,LN
           INCRI=INCRM+I
380         B(INCRI) = B(INCRI) - ROWMLT*BOTBLK(I,J)
390 CONTINUE
400 CONTINUE
   RETURN
   END
   SUBROUTINE CRSLVE(TOPBLK,ARRAY,BOTBLK,IPIVOT,B)
   IMPLICIT REAL*8(A-H,O-Z)
   PARAMETER(LM=10,LN=1,KK=10,
  ,  LM1=LM+1,LN1=LN+1,KK1=KK+1,
  ,  LC=LM+LN,LC1=LC+1,LA=2*LC,N=LC*KK1)
   DIMENSION TOPBLK(LM,LC),ARRAY(LC,LA,KK),BOTBLK(LN,LC)
   DIMENSION IPIVOT(N),B(N)
   LB=LC+LM
   LB1=LB+1
   INCR=KK*LC
   INCRM=INCR+LM
   DO 210 L=1,LN
     J = LC1 - L
     INCRJ = INCR + J
     LN1L = LN1 - L
     B(INCRJ) = B(INCRJ)/BOTBLK(LN1L,J)
     IF (L.EQ.LN) GO TO 200
     BINCRJ = B(INCRJ)
     LNL = LN - L
     DO 190 I=1,LNL
       INCRI = INCRM + I
       B(INCRI) = B(INCRI) - BOTBLK(I,J)*BINCRJ
190   CONTINUE
```

```
 200   CONTINUE
 210 CONTINUE
     DO 300 LK=1,KK
       K = KK1 - LK
       INCR = INCR - LC
       DO 225 L1=LN1,LC
         I = LC + LN1 - L1
         IM = I + LM
         INCRN = INCR + IM
         DOTPRD = B(INCRN)
         DO 220 J=LB1,LA
           INCRJ = INCR + J
           DOTPRD = DOTPRD - ARRAY(I,J,K)*B(INCRJ)
 220     CONTINUE
         B(INCRN) = DOTPRD
 225 CONTINUE
       DO 240 L1=LN1,LC
         I = LC + LN1 - L1
         IM = I + LM
         INCRN = INCR + IM
         IPVTN = IPIVOT(INCRN)
         IF (INCRN.EQ.IPVTN) GO TO 230
         SWAP = B(INCRN)
         B(INCRN) = B(IPVTN)
         B(IPVTN) = SWAP
 230     CONTINUE
 240   CONTINUE
       INCRM = INCR + LM
       DO 260 J=LC1,LA
         INCRJ = INCR + J
         BINCRJ = B(INCRJ)
         DO 250 I=1,LN
           INCRI = INCRM + I
           B(INCRI) = B(INCRI) - ARRAY(I,J,K)*BINCRJ
 250     CONTINUE
 260   CONTINUE
       DO 290 L=1,LN
         J = LC1 - L
         INCRJ = INCR + J
         LN1L = LN1 - L
         B(INCRJ) = B(INCRJ)/ARRAY(LN1L,J,K)
         IF (L.EQ.LN) GO TO 280
         BINCRJ = B(INCRJ)
         LNL = LN - L
         DO 270 I=1,LNL
           INCRI = INCRM + I
           B(INCRI) = B(INCRI) - ARRAY(I,J,K)*BINCRJ
 270     CONTINUE
 280     CONTINUE
 290   CONTINUE
 300 CONTINUE
     DO 315 L=1,LM
       I = LM1 - L
       DOTPRD = B(I)
       DO 310 J=LM1,LC
         DOTPRD = DOTPRD - TOPBLK(I,J)*B(J)
 310   CONTINUE
       B(I) = DOTPRD
 315 CONTINUE
     DO 330 L=1,LM
       I = LM1 - L
       IPVTI = IPIVOT(I)
       IF (I.EQ.IPVTI) GO TO 320
       SWAP = B(I)
       B(I) = B(IPVTI)
       B(IPVTI) = SWAP
 320   CONTINUE
 330 CONTINUE
     RETURN
     END
```

## BCBR.FOR Fortran program

```
CCCCC
```

```
C    BCBR
CCCCC
      IMPLICIT REAL*8(A-H,O-Z)
      PARAMETER(LM=10,LN=1,KK=10,MM=1000000,
     ,  LM1=LM+1,LN1=LN+1,KK1=KK+1,
     ,  LC=LM+LN,LC1=LC+1,LA=2*LC,LA1=LA+1,N=LC*KK1)
      DIMENSION TOPBLK(LM,LC),ARRAY(LC,LA,KK),BOTBLK(LN,LC)
      DIMENSION IPIVOT(N),B(N)
      DIMENSION A1(LM,LC1),A2(LN,LC1),A3(LC,LA1)
      OPEN(1,FILE='TESTA')
      READ(1,1000) ((A1(I,J),J=1,LC1),I=1,LM)
      READ(1,1000) ((A2(I,J),J=1,LC1),I=1,LN)
      CLOSE(1)
      OPEN(2,FILE='TESTB')
      READ(2,2000) ((A3(I,J),J=1,LA1),I=1,LC)
      CLOSE(2)
      DO 4 M=1,MM
      DO 1 I=1,LM
      B(I)=A1(I,LC1)
      DO 1 J=1,LC
    1 TOPBLK(I,J)=A1(I,J)
      DO 2 K=1,KK
      DO 2 I=1,LC
      B(I+K*LC-LN)=A3(I,LA1)
      DO 2 J=1,LA
    2 ARRAY(I,J,K)=A3(I,J)
      DO 3 I=1,LN
      B(I+N-LN)=A2(I,LC1)
      DO 3 J=1,LC
    3 BOTBLK(I,J)=A2(I,J)
      CALL CRDCMP(TOPBLK,ARRAY,BOTBLK,IPIVOT,B)
      CALL CRSLVE(TOPBLK,ARRAY,BOTBLK,IPIVOT,B)
    4 CONTINUE
      OPEN(3,FILE='TESTC')
      WRITE (3,3000) (B(I),I=1,N)
      CLOSE(3)
      STOP
 1000 FORMAT(5X/(11F6.0,F8.0))
 2000 FORMAT(5X/(22F6.0,F8.0))
 3000 FORMAT(3E25.15)
      END
      SUBROUTINE CRDCMP(TOPBLK,ARRAY,BOTBLK,IPIVOT,B)
      IMPLICIT REAL*8(A-H,O-Z)
      PARAMETER(LM=10,LN=1,KK=10,
     ,  LM1=LM+1,LN1=LN+1,KK1=KK+1,
     ,  LC=LM+LN,LC1=LC+1,LA=2*LC,N=LC*KK1)
      DIMENSION TOPBLK(LM,LC),ARRAY(LC,LA,KK),BOTBLK(LN,LC)
      DIMENSION IPIVOT(N),B(N)
      LB=LC+LM
      LB1=LB+1
      DO 80 I=1,LM
        I1 = I + 1
        IPVT = I
        COLMAX = DABS(TOPBLK(I,I))
        DO 30 J=I1,LC
          TEMPIV = DABS(TOPBLK(I,J))
          IF (TEMPIV.LE.COLMAX) GO TO 30
          IPVT = J
          COLMAX = TEMPIV
   30   CONTINUE
        IPIVOT(I) = IPVT
        IF (IPVT.EQ.I) GO TO 60
        DO 40 L=1,LM
          SWAP = TOPBLK(L,IPVT)
          TOPBLK(L,IPVT) = TOPBLK(L,I)
          TOPBLK(L,I) = SWAP
   40   CONTINUE
        DO 50 L=1,LC
          SWAP = ARRAY(L,IPVT,1)
          ARRAY(L,IPVT,1) = ARRAY(L,I,1)
          ARRAY(L,I,1) = SWAP
   50   CONTINUE
   60   CONTINUE
        CII = TOPBLK(I,I)
        DO 70 J=I1,LC
```

```
        CIJ = TOPBLK(I,J)/CII
        TOPBLK(I,J) = CIJ
        DO 70 L=I1,LM
 70       TOPBLK(L,J) = TOPBLK(L,J) - TOPBLK(L,I)*CIJ
    BI=B(I)/CII
    B(I)=BI
    DO 75 L=I1,LM
 75 B(L)=B(L)-TOPBLK(L,I)*BI
 80 CONTINUE
    DO 110 L=LM,1,-1
    L1=L-1
    BL=B(L)
    DO 90 I=1,L1
    CIL=TOPBLK(I,L)
    DO 85 J=LM1,LC
 85 TOPBLK(I,J)=TOPBLK(I,J)-CIL*TOPBLK(L,J)
 90 B(I)=B(I)-CIL*BL
    DO 100 I=1,LC
    IM=I+LM
    CIL=ARRAY(I,L,1)
    DO 95 J=LM1,LC
 95 ARRAY(I,J,1) = ARRAY(I,J,1) - CIL*TOPBLK(L,J)
100 B(IM)=B(IM)-CIL*BL
110 CONTINUE
    INCR = 0
    DO 320 K=1,KK
      K1 = K + 1
        INCRM=INCR+LM
        INCRC=INCR+LC
        INCRB=INCR+LB
      DO 170 J=LM1,LC
        J1 = J + 1
        JN = J - LM
        JN1 = JN + 1
        INCRJ=INCR+J
        IPVT = JN
        ROWMAX = DABS(ARRAY(JN,J,K))
        DO 120 I=JN1,LC
          TEMPIV = DABS(ARRAY(I,J,K))
          IF (TEMPIV.LE.ROWMAX) GO TO 120
          IPVT = I
          ROWMAX = TEMPIV
120     CONTINUE
        INCRP=INCRM+IPVT
        IPIVOT(INCRJ)=INCRP
        IF (IPVT.EQ.JN) GO TO 140
        DO 130 L=LM1,LA
          SWAP = ARRAY(IPVT,L,K)
          ARRAY(IPVT,L,K) = ARRAY(JN,L,K)
          ARRAY(JN,L,K) = SWAP
130     CONTINUE
          SWAP = B(INCRP)
          B(INCRP) = B(INCRJ)
          B(INCRJ) = SWAP
140     CONTINUE
        ROWPIV = ARRAY(JN,J,K)
        DO 150 I=JN1,LC
          ARRAY(I,J,K) = ARRAY(I,J,K)/ROWPIV
150     CONTINUE
        DO 160 L=J1,LC
          ROWMLT = ARRAY(JN,L,K)
          DO 160 I=JN1,LC
160         ARRAY(I,L,K) = ARRAY(I,L,K) - ROWMLT*ARRAY(I,J,K)
170   CONTINUE
    L2=LC
    DO 180 L=LN,2,-1
    L1=L2-1
    DO 175 J=L1,1,-1
    ANJ=ARRAY(L,J,K)
    DO 175 I=LN1,LC
175 ARRAY(I,J,K)=ARRAY(I,J,K)-ARRAY(I,L2,K)*ANJ
180 L2=L1
    DO 182 L=1,LN
    LLM=L+LM
    BL=B(L+INCRM)
```

```
   DO 182 I=LN1,LC
   AIN=ARRAY(I,LLM,K)
   DO 181 J=LC1,LA
181 ARRAY(I,J,K)=ARRAY(I,J,K)-AIN*ARRAY(L,J,K)
   IINCRM=I+INCRM
182 B(IINCRM)=B(IINCRM)-AIN*BL
     DO 260 I=LN1,LC
       IM = I + LM
       IM1 = IM + 1
       IMC = IM - LC
       INCRN = INCR + IM
       IN=I-LN
       I1 = I + 1
       IPVT = IM
       COLMAX = DABS(ARRAY(I,IPVT,K))
       DO 190 J=IM1,LA
         TEMPIV = DABS(ARRAY(I,J,K))
         IF (TEMPIV.LE.COLMAX) GO TO 190
         IPVT = J
         COLMAX = TEMPIV
190     CONTINUE
       IPIVOT(INCRN) = INCR + IPVT
       IF (IPVT.EQ.IM) GO TO 240
       DO 200 L=LN1,LC
         SWAP = ARRAY(L,IPVT,K)
         ARRAY(L,IPVT,K) = ARRAY(L,IM,K)
         ARRAY(L,IM,K) = SWAP
200     CONTINUE
       IPC = IPVT - LC
       IF (K.EQ.KK) GO TO 220
       DO 210 L=1,LC
         SWAP = ARRAY(L,IPC,K1)
         ARRAY(L,IPC,K1) = ARRAY(L,IMC,K1)
         ARRAY(L,IMC,K1) = SWAP
210     CONTINUE
       GO TO 240
220     CONTINUE
       DO 230 L=1,LN
         SWAP = BOTBLK(L,IPC)
         BOTBLK(L,IPC) = BOTBLK(L,IMC)
         BOTBLK(L,IMC) = SWAP
230     CONTINUE
240     CONTINUE
       CII = ARRAY(I,IM,K)
       DO 250 J=IM1,LA
         CIJ = ARRAY(I,J,K)/CII
         ARRAY(I,J,K) = CIJ
         DO 250 L=I1,LC
250          ARRAY(L,J,K) = ARRAY(L,J,K) - ARRAY(L,IM,K)*CIJ
   BI=B(INCRN)/CII
   B(INCRN)=BI
   DO 255 L=I1,LC
   LINCRM=L+INCRM
255 B(LINCRM)=B(LINCRM)-ARRAY(L,IM,K)*BI
260 CONTINUE
   DO 300 L=LM,1,-1
   L1=L-1
   L1N=L1+LN
   BL=B(L+INCRC)
   DO 270 I=LN1,L1N
   IINCRM=I+INCRM
   CIL=ARRAY(I,L+LC,K)
   DO 265 J=LB1,LA
265 ARRAY(I,J,K)=ARRAY(I,J,K)-CIL*ARRAY(L+LN,J,K)
270 B(IINCRM)=B(IINCRM)-CIL*BL
   IF(K.LT.KK) THEN
   DO 280 I=1,LC
   IINCRB=I+INCRB
   CIL=ARRAY(I,L,K1)
   DO 275 J=LM1,LC
275 ARRAY(I,J,K1) = ARRAY(I,J,K1) - CIL*ARRAY(L+LN,J+LC,K)
280 B(IINCRB)=B(IINCRB)-CIL*BL
   ELSE
   DO 290 I=1,LN
   IINCRB=I+INCRB
```

```
   CIL=BOTBLK(I,L)
   DO 285 J=LM1,LC
285 BOTBLK(I,J) = BOTBLK(I,J) - CIL*ARRAY(L+LN,J+LC,K)
290 B(IINCRB)=B(IINCRB)-CIL*BL
   ENDIF
300 CONTINUE
     INCR = INCRC
320 CONTINUE
   IF (LN.EQ.1) GO TO 400
   INCRM=INCR+LM
   LC2=LC-1
   DO 390 J=LM1,LC2
     J1 = J + 1
     JN = J - LM
     JN1 = JN + 1
     INCRJ=INCR+J
     IPVT = JN
     ROWMAX = DABS(BOTBLK(JN,J))
     DO 330 I=JN1,LN
       TEMPIV = DABS(BOTBLK(I,J))
       IF (TEMPIV.LE.ROWMAX) GO TO 330
       IPVT = I
       ROWMAX = TEMPIV
330   CONTINUE
     INCRP=INCRM+IPVT
     IPIVOT(INCRJ)=INCRP
     IF (IPVT.EQ.JN) GO TO 350
     DO 340 L=LM1,LC
       SWAP = BOTBLK(IPVT,L)
       BOTBLK(IPVT,L) = BOTBLK(JN,L)
       BOTBLK(JN,L) = SWAP
340   CONTINUE
         SWAP = B(INCRP)
         B(INCRP) = B(INCRJ)
         B(INCRJ) = SWAP
350   CONTINUE
     ROWPIV = BOTBLK(JN,J)
     DO 360 I=JN1,LN
360     BOTBLK(I,J) = BOTBLK(I,J)/ROWPIV
     DO 370 L=J1,LC
       ROWMLT = BOTBLK(JN,L)
       DO 370 I=JN1,LN
370       BOTBLK(I,L) = BOTBLK(I,L) - ROWMLT*BOTBLK(I,J)
         ROWMLT = B(INCRJ)
         DO 380 I=JN1,LN
           INCRI=INCRM+I
380         B(INCRI) = B(INCRI) - ROWMLT*BOTBLK(I,J)
390 CONTINUE
400 CONTINUE
   RETURN
   END
   SUBROUTINE CRSLVE(TOPBLK,ARRAY,BOTBLK,IPIVOT,B)
   IMPLICIT REAL*8(A-H,O-Z)
   PARAMETER(LM=10,LN=1,KK=10,
  ,  LM1=LM+1,LN1=LN+1,KK1=KK+1,
  ,  LC=LM+LN,LC1=LC+1,LA=2*LC,N=LC*KK1)
   DIMENSION TOPBLK(LM,LC),ARRAY(LC,LA,KK),BOTBLK(LN,LC)
   DIMENSION IPIVOT(N),B(N)
   LB=LC+LM
   LB1=LB+1
   INCR=KK*LC
   INCRM=INCR+LM
   DO 210 L=1,LN
     J = LC1 - L
     INCRJ = INCR + J
     LN1L = LN1 - L
     B(INCRJ) = B(INCRJ)/BOTBLK(LN1L,J)
     IF (L.EQ.LN) GO TO 200
     BINCRJ = B(INCRJ)
     LNL = LN - L
     DO 190 I=1,LNL
       INCRI = INCRM + I
       B(INCRI) = B(INCRI) - BOTBLK(I,J)*BINCRJ
190   CONTINUE
200   CONTINUE
```

```
210 CONTINUE
    DO 300 LK=1,KK
      K = KK1 - LK
      INCR = INCR - LC
      DO 225 L1=LN1,LC
        I = LC + LN1 - L1
        IM = I + LM
        INCRN = INCR + IM
        DOTPRD = B(INCRN)
        DO 220 J=LB1,LA
          INCRJ = INCR + J
          DOTPRD = DOTPRD - ARRAY(I,J,K)*B(INCRJ)
220     CONTINUE
        B(INCRN) = DOTPRD
225 CONTINUE
      DO 240 L1=LN1,LC
        I = LC + LN1 - L1
        IM = I + LM
        INCRN = INCR + IM
        IPVTN = IPIVOT(INCRN)
        IF (INCRN.EQ.IPVTN) GO TO 230
        SWAP = B(INCRN)
        B(INCRN) = B(IPVTN)
        B(IPVTN) = SWAP
230     CONTINUE
240   CONTINUE
      INCRM = INCR + LM
      DO 260 J=LC1,LA
        INCRJ = INCR + J
        BINCRJ = B(INCRJ)
        DO 250 I=1,LN
          INCRI = INCRM + I
          B(INCRI) = B(INCRI) - ARRAY(I,J,K)*BINCRJ
250     CONTINUE
260   CONTINUE
    DO 295 L=2,LN
    L1=L-1
    L1LM=L1+LM
    BINCRL=B(INCRM+L1)
    DO 295 I=L,LN
    INCRI=INCRM+I
295 B(INCRI)=B(INCRI)-BINCRL*ARRAY(I,L1LM,K)
      DO 290 L=1,LN
        J = LC1 - L
        INCRJ = INCR + J
        LN1L = LN1 - L
        B(INCRJ) = B(INCRJ)/ARRAY(LN1L,J,K)
        IF (L.EQ.LN) GO TO 280
        BINCRJ = B(INCRJ)
        LNL = LN - L
        DO 270 I=1,LNL
          INCRI = INCRM + I
          B(INCRI) = B(INCRI) - ARRAY(I,J,K)*BINCRJ
270     CONTINUE
280     CONTINUE
290   CONTINUE
300 CONTINUE
    DO 315 L=1,LM
      I = LM1 - L
      DOTPRD = B(I)
      DO 310 J=LM1,LC
        DOTPRD = DOTPRD - TOPBLK(I,J)*B(J)
310   CONTINUE
      B(I) = DOTPRD
315 CONTINUE
    DO 330 L=1,LM
      I = LM1 - L
      IPVTI = IPIVOT(I)
      IF (I.EQ.IPVTI) GO TO 320
      SWAP = B(I)
      B(I) = B(IPVTI)
      B(IPVTI) = SWAP
320   CONTINUE
330 CONTINUE
    RETURN
```

```
      END
```

## DBTC.FOR Fortran program

```
CCCCC
C    DBTC
CCCCC
      IMPLICIT REAL*8(A-H,O-Z)
      PARAMETER(LM=10,LN=1,KK=10,MM=1000000,
     ,  LM1=LM+1,LN1=LN+1,KK1=KK+1,
     ,  LC=LM+LN,LC1=LC+1,LA=2*LC,LA1=LA+1,N=LC*KK1)
      DIMENSION TOPBLK(LM,LC),ARRAY(LC,LA,KK),BOTBLK(LN,LC)
      DIMENSION IPIVOT(N),B(N)
      DIMENSION A1(LM,LC1),A2(LN,LC1),A3(LC,LA1)
      OPEN(1,FILE='TESTA')
      READ(1,1000) ((A1(I,J),J=1,LC1),I=1,LM)
      READ(1,1000) ((A2(I,J),J=1,LC1),I=1,LN)
      CLOSE(1)
      OPEN(2,FILE='TESTB')
      READ(2,2000) ((A3(I,J),J=1,LA1),I=1,LC)
      CLOSE(2)
      DO 4 M=1,MM
      DO 1 I=1,LM
      B(I)=A1(I,LC1)
      DO 1 J=1,LC
    1 TOPBLK(I,J)=A1(I,J)
      DO 2 K=1,KK
      DO 2 I=1,LC
      B(I+K*LC-LN)=A3(I,LA1)
      DO 2 J=1,LA
    2 ARRAY(I,J,K)=A3(I,J)
      DO 3 I=1,LN
      B(I+N-LN)=A2(I,LC1)
      DO 3 J=1,LC
    3 BOTBLK(I,J)=A2(I,J)
      CALL CRDCMP(TOPBLK,ARRAY,BOTBLK,IPIVOT,B)
      CALL CRSLVE(TOPBLK,ARRAY,BOTBLK,IPIVOT,B)
    4 CONTINUE
      OPEN(3,FILE='TESTC')
      WRITE (3,3000) (B(I),I=1,N)
      CLOSE(3)
      STOP
 1000 FORMAT(5X/(11F6.0,F8.0))
 2000 FORMAT(5X/(22F6.0,F8.0))
 3000 FORMAT(3E25.15)
      END
      SUBROUTINE CRDCMP(TOPBLK,ARRAY,BOTBLK,IPIVOT,B)
      IMPLICIT REAL*8(A-H,O-Z)
      PARAMETER(LM=10,LN=1,KK=10,
     ,  LM1=LM+1,LN1=LN+1,KK1=KK+1,
     ,  LC=LM+LN,LC1=LC+1,LA=2*LC,N=LC*KK1)
      DIMENSION TOPBLK(LM,LC),ARRAY(LC,LA,KK),BOTBLK(LN,LC)
      DIMENSION IPIVOT(N),B(N)
      LB=LC+LM
      LB1=LB+1
      DO 80 I=1,LM
        I1 = I + 1
        IPVT = I
        COLMAX = DABS(TOPBLK(I,I))
        DO 30 J=I1,LC
          TEMPIV = DABS(TOPBLK(I,J))
          IF (TEMPIV.LE.COLMAX) GO TO 30
          IPVT = J
          COLMAX = TEMPIV
   30   CONTINUE
        IPIVOT(I) = IPVT
        IF (IPVT.EQ.I) GO TO 60
        DO 40 L=1,LM
          SWAP = TOPBLK(L,IPVT)
          TOPBLK(L,IPVT) = TOPBLK(L,I)
          TOPBLK(L,I) = SWAP
   40   CONTINUE
        DO 50 L=1,LC
          SWAP = ARRAY(L,IPVT,1)
```

```
       ARRAY(L,IPVT,1) = ARRAY(L,I,1)
       ARRAY(L,I,1) = SWAP
 50   CONTINUE
 60   CONTINUE
      CII = TOPBLK(I,I)
      DO 70 J=I1,LC
       CIJ = TOPBLK(I,J)/CII
       TOPBLK(I,J) = CIJ
       DO 70 L=I1,LM
 70       TOPBLK(L,J) = TOPBLK(L,J) - TOPBLK(L,I)*CIJ
    BI=B(I)/CII
    B(I)=BI
    DO 75 L=I1,LM
 75 B(L)=B(L)-TOPBLK(L,I)*BI
 80 CONTINUE
    DO 110 L=LM,1,-1
    L1=L-1
    BL=B(L)
    DO 90 I=1,L1
    CIL=TOPBLK(I,L)
    DO 85 J=LM1,LC
 85 TOPBLK(I,J)=TOPBLK(I,J)-CIL*TOPBLK(L,J)
 90 B(I)=B(I)-CIL*BL
    DO 100 I=1,LC
    IM=I+LM
    CIL=ARRAY(I,L,1)
    DO 95 J=LM1,LC
 95 ARRAY(I,J,1) = ARRAY(I,J,1) - CIL*TOPBLK(L,J)
100 B(IM)=B(IM)-CIL*BL
110 CONTINUE
    INCR = 0
    DO 320 K=1,KK
      K1 = K + 1
       INCRM=INCR+LM
       INCRC=INCR+LC
       INCRB=INCR+LB
      DO 180 J=LM1,LC
       J1 = J + 1
       JN = J - LM
       JN1 = JN + 1
       INCRJ=INCR+J
       IPVT = JN
       ROWMAX = DABS(ARRAY(JN,J,K))
       DO 120 I=JN1,LC
         TEMPIV = DABS(ARRAY(I,J,K))
         IF (TEMPIV.LE.ROWMAX) GO TO 120
         IPVT = I
         ROWMAX = TEMPIV
120     CONTINUE
       INCRP=INCRM+IPVT
       IPIVOT(INCRJ)=INCRP
       IF (IPVT.EQ.JN) GO TO 140
       DO 130 L=J,LA
         SWAP = ARRAY(IPVT,L,K)
         ARRAY(IPVT,L,K) = ARRAY(JN,L,K)
         ARRAY(JN,L,K) = SWAP
130     CONTINUE
         SWAP = B(INCRP)
         B(INCRP) = B(INCRJ)
         B(INCRJ) = SWAP
140     CONTINUE
       ROWPIV = ARRAY(JN,J,K)
       DO 150 I=JN1,LC
         ARRAY(I,J,K) = ARRAY(I,J,K)/ROWPIV
150     CONTINUE
       DO 160 L=J1,LA
         ROWMLT = ARRAY(JN,L,K)
         DO 160 I=JN1,LC
160         ARRAY(I,L,K) = ARRAY(I,L,K) - ROWMLT*ARRAY(I,J,K)
         ROWMLT = B(INCRJ)
         DO 170 I=JN1,LC
           INCRI=INCRM+I
170         B(INCRI) = B(INCRI) - ROWMLT*ARRAY(I,J,K)
180   CONTINUE
      DO 260 I=LN1,LC
```

```
       IM = I + LM
       IM1 = IM + 1
       IMC = IM - LC
       INCRN = INCR + IM
       IN=I-LN
       I1 = I + 1
       IPVT = IM
       COLMAX = DABS(ARRAY(I,IPVT,K))
       DO 190 J=IM1,LA
         TEMPIV = DABS(ARRAY(I,J,K))
         IF (TEMPIV.LE.COLMAX) GO TO 190
         IPVT = J
         COLMAX = TEMPIV
 190     CONTINUE
       IPIVOT(INCRN) = INCR + IPVT
       IF (IPVT.EQ.IM) GO TO 240
       DO 200 L=1,LC                                    CC
         SWAP = ARRAY(L,IPVT,K)
         ARRAY(L,IPVT,K) = ARRAY(L,IM,K)
         ARRAY(L,IM,K) = SWAP
 200     CONTINUE
       IPC = IPVT - LC
       IF (K.EQ.KK) GO TO 220
       DO 210 L=1,LC
         SWAP = ARRAY(L,IPC,K1)
         ARRAY(L,IPC,K1) = ARRAY(L,IMC,K1)
         ARRAY(L,IMC,K1) = SWAP
 210     CONTINUE
       GO TO 240
 220     CONTINUE
       DO 230 L=1,LN
         SWAP = BOTBLK(L,IPC)
         BOTBLK(L,IPC) = BOTBLK(L,IMC)
         BOTBLK(L,IMC) = SWAP
 230     CONTINUE
 240     CONTINUE
       CII = ARRAY(I,IM,K)
       DO 254 J=IM1,LA
         CIJ = ARRAY(I,J,K)/CII
         ARRAY(I,J,K) = CIJ
         DO 250 L=I1,LC
 250          ARRAY(L,J,K) = ARRAY(L,J,K) - ARRAY(L,IM,K)*CIJ
 254       CONTINUE
     BI=B(INCRN)/CII
     B(INCRN)=BI
     DO 255 L=I1,LC
     LINCRM=L+INCRM
 255 B(LINCRM)=B(LINCRM)-ARRAY(L,IM,K)*BI
 260 CONTINUE
     DO 270 L=LM,1,-1
     L1=L-1
     L1N=L1+LN
     BL=B(L+INCRC)
     DO 270 I=LN1,L1N
     IINCRM=I+INCRM
     CIL=ARRAY(I,L+LC,K)
     DO 265 J=LB1,LA
 265 ARRAY(I,J,K)=ARRAY(I,J,K)-CIL*ARRAY(L+LN,J,K)
 270 B(IINCRM)=B(IINCRM)-CIL*BL
       DO 259 I=LN1,LC
       IM=I+LM
       IN=I-LN
       DO 2540 J=LB1,LA
         CIJ = ARRAY(I,J,K)
         DO 251 L=1,LN
 251          ARRAY(L,J,K) = ARRAY(L,J,K) - ARRAY(L,IM,K)*CIJ
         JC = J - LC
         IF (K.LT.KK) THEN
         DO 252 L=1,LC
 252        ARRAY(L,JC,K1) = ARRAY(L,JC,K1) -ARRAY(L,IN,K1)*CIJ
         ELSE
         DO 253 L=1,LN
 253        BOTBLK(L,JC) = BOTBLK(L,JC) -BOTBLK(L,IN)*CIJ
         ENDIF
2540 CONTINUE
```

```
   BI=B(I+INCRM)
   DO 256 L=1,LN
   LINCRM=L+INCRM
256 B(LINCRM)=B(LINCRM)-ARRAY(L,IM,K)*BI
   IF (K.LT.KK) THEN
   DO 257 L=1,LC
   LINCRB=L+INCRB
257 B(LINCRB)=B(LINCRB)-ARRAY(L,IN,K1)*BI
   ELSE
   DO 258 L=1,LN
   LINCRB=L+INCRB
258 B(LINCRB)=B(LINCRB)-BOTBLK(L,IN)*BI
   ENDIF
259       CONTINUE
   DO 310 L=1,LN
     LL = LC1 - L
     INCRL = INCR + LL
     LN1L = LN1 - L
     LNL = LN - L
     ALL=ARRAY(LN1L,LL,K)
     BINCRL = B(INCRL)/ALL
     B(INCRL)=BINCRL
     DO 301  J=LB1,LA
     ALJ=ARRAY(LN1L,J,K)/ALL
     ARRAY(LN1L,J,K)=ALJ
     DO 301 I=1,LNL
301   ARRAY(I,J,K)=ARRAY(I,J,K)-ARRAY(I,LL,K)*ALJ
     DO 302 I=1,LNL
       INCRI = INCRM + I
302     B(INCRI) = B(INCRI) - ARRAY(I,LL,K)*BINCRL
310 CONTINUE
     INCR = INCRC
320 CONTINUE
   IF (LN.EQ.1) GO TO 400
   INCRM=INCR+LM
   LC2=LC-1
   DO 390 J=LM1,LC2
     J1 = J + 1
     JN2 = J - LM
     JN21 = JN2 + 1
     INCRJ=INCR+J
     IPVT = JN2
     ROWMAX = DABS(BOTBLK(JN2,J))
     DO 330 I=JN21,LN
       TEMPIV = DABS(BOTBLK(I,J))
       IF (TEMPIV.LE.ROWMAX) GO TO 330
       IPVT = I
       ROWMAX = TEMPIV
330   CONTINUE
     INCRP=INCRM+IPVT
     IPIVOT(INCRJ)=INCRP
     IF (IPVT.EQ.JN2) GO TO 350
     DO 340 L=J,LC
       SWAP = BOTBLK(IPVT,L)
       BOTBLK(IPVT,L) = BOTBLK(JN2,L)
       BOTBLK(JN2,L) = SWAP
340   CONTINUE
         SWAP = B(INCRP)
         B(INCRP) = B(INCRJ)
         B(INCRJ) = SWAP
350   CONTINUE
     ROWPIV = BOTBLK(JN2,J)
     DO 360 I=JN21,LN
360     BOTBLK(I,J) = BOTBLK(I,J)/ROWPIV
     DO 370 L=J1,LC
       ROWMLT = BOTBLK(JN2,L)
       DO 370 I=JN21,LN
370       BOTBLK(I,L) = BOTBLK(I,L) - ROWMLT*BOTBLK(I,J)
         ROWMLT = B(INCRJ)
         DO 380 I=JN21,LN
           INCRI=INCRM+I
380         B(INCRI) = B(INCRI) - ROWMLT*BOTBLK(I,J)
390 CONTINUE
400 CONTINUE
   RETURN
```

```
      END
      SUBROUTINE CRSLVE(TOPBLK,ARRAY,BOTBLK,IPIVOT,B)
      IMPLICIT REAL*8(A-H,O-Z)
      PARAMETER(LM=10,LN=1,KK=10,
     ,  LM1=LM+1,LN1=LN+1,KK1=KK+1,
     ,  LC=LM+LN,LC1=LC+1,LA=2*LC,N=LC*KK1)
      DIMENSION TOPBLK(LM,LC),ARRAY(LC,LA,KK),BOTBLK(LN,LC)
      DIMENSION IPIVOT(N),B(N)
      LB=LC+LM
      LB1=LB+1
      INCR=KK*LC
      INCRM=INCR+LM
      DO 210 L=1,LN
        J = LC1 - L
        INCRJ = INCR + J
        LN1L = LN1 - L
        B(INCRJ) = B(INCRJ)/BOTBLK(LN1L,J)
        IF (L.EQ.LN) GO TO 200
        BINCRJ = B(INCRJ)
        LNL = LN - L
        DO 190 I=1,LNL
          INCRI = INCRM + I
          B(INCRI) = B(INCRI) - BOTBLK(I,J)*BINCRJ
190   CONTINUE
200   CONTINUE
210 CONTINUE
      DO 300 LK=1,KK
        K = KK1 - LK
        INCR = INCR - LC
        DO 225 L=1,LC
          IC = LC1-L
          IB = LB1 - L
          INCRI = INCR + IB
          DOTPRD = B(INCRI)
          DO 220 J=LB1,LA
            INCRJ = INCR + J
            DOTPRD = DOTPRD - ARRAY(IC,J,K)*B(INCRJ)
220     CONTINUE
          B(INCRI) = DOTPRD
225 CONTINUE
        DO 240 L=1,LM
          IB = LB1-L
          INCRI = INCR + IB
          IPVTI = IPIVOT(INCRI)
          IF (INCRI.EQ.IPVTI) GO TO 230
          SWAP = B(INCRI)
          B(INCRI) = B(IPVTI)
          B(IPVTI) = SWAP
230     CONTINUE
240   CONTINUE
300 CONTINUE
      DO 315 L=1,LM
        I = LM1 - L
        DOTPRD = B(I)
        DO 310 J=LM1,LC
          DOTPRD = DOTPRD - TOPBLK(I,J)*B(J)
310   CONTINUE
        B(I) = DOTPRD
315 CONTINUE
      DO 330 L=1,LM
        I = LM1 - L
        IPVTI = IPIVOT(I)
        IF (I.EQ.IPVTI) GO TO 320
        SWAP = B(I)
        B(I) = B(IPVTI)
        B(IPVTI) = SWAP
320   CONTINUE
330 CONTINUE
      RETURN
      END
```

## TESTA input data

```
TOP BLOCK
  .10 -.22 -.24 -.42 -.37 -.77 -.99 -.96 -.89  .85 -.28 -6.412
 -.63  .09 -.10 -.07  .51 -.02  .01 -.52  .07  .48 -.54 -0.968
  .32 -.29  .02 -.81  .29  .00 -.05 -.91  .00  .00  .69 -0.869
 -.25 -.09 -.91 -.17 -.46 -.92 -.14  .98 -.34  .70 -.53 -2.586
  .76 -.90 -.64 -.08  .95  .15  .15 -.46 -.48  .93 -.39  0.034
 -.06 -.72 -.91 -.14  .36 -.69 -.40 -.93 -.61 -.97 -.12 -8.059
 -.21  .54 -.53  .97  .91  .58 -.32  .27  .33  .72 -.20  4.662
 -.57  .04  .26 -.04  .69 -.65 -.57  .83 -.42 -.56 -.18 -1.956
  .89 -.62 -.07 -.63  .28 -.54 -.29  .52  .67  .00 -.68 -0.847
  .10 -.01 -.25 -.22  .06  .81  .11  .56  .05  .63 -.43  2.357
BOTTOM BLOCK
  .88  .48  .52 -.87  .71  .51  .52 -.33 -.46 -.33  .85  3.176
```

## TESTB input data

```
RECURRENT BLOCK 11X(22+1)
  .22 -.05  .87  .28  .04  .68  .39  .25 -.64 -.87 -.62  .95  .29 -.73 -.27 -.90  .18  .94  .35 -.33 -.88  .39 -0.682
 -.06 -.40 -.83 -.33  .31 -.93  .20  .02 -.85  .97  .61  .16 -.42 -.69 -.07  .10 -.53  .33  .03 -.92  .85 -.08 -2.497
  .51  .60 -.94 -.58 -.09 -.14 -.74  .24 -.87 -.07  .96  .26  .66  .26 -.94 -.77 -.56  .55  .88 -.12 -.30 -.49 -2.835
  .49 -.78  .81  .64 -.82  .46  .67 -.07 -.29 -.31 -.25 -.70 -.38  .81 -.30 -.76  .07 -.06 -.27  .98  .18  .17  0.716
  .53 -.70  .49 -.88  .48  .77  .77 -.89  .31  .23  .42 -.09  .47 -.13 -.58 -.19  .24 -.46  .84  .44 -.26  .42  4.026
 -.86 -.18 -.67  .30  .04  .20 -.02  .84  .39  .01  .34  .23 -.68 -.58  .65  .14  .61 -.10 -.91  .91 -.89  .64  1.943
 -.16 -.91  .53  .31 -.20 -.18 -.59 -.79  .69  .33  .52 -.13 -.16  .19 -.04 -.14  .06 -.30 -.25  .38  .00  .92  1.333
  .82  .20  .40  .44 -.25 -.35  .88 -.27 -.48 -.18 -.86 -.59  .51  .82 -.47  .92  .17  .53 -.82 -.25  .38  .24  1.333
  .34 -.04 -.21 -.69 -.27 -.15  .39  .60  .18 -.71 -.94 -.57  .38  .56  .18 -.36  .67 -.47 -.60 -.55 -.18 -.83 -6.226
 -.76  .90  .76 -.94  .45 -.60 -.66 -.89  .32 -.37 -.39  .74  .76 -.26 -.18  .28  .29 -.06  .20 -.09  .56 -.66 -2.143
 -.87 -.94 -.11 -.22  .50 -.59  .81  .76 -.59 -.14  .53  .24 -.53 -.81  .70 -.18  .56 -.84 -.62  .05  .72  .17 -0.604
```